\catcode`\@=11
\hsize 15 truecm \hoffset 0.46 truecm
\vsize 23 truecm \voffset -0.04 truecm
\baselineskip=14pt plus 1pt
\parskip=2pt plus 1pt minus 1pt
\parindent=0mm  \hfuzz=30pt
\pretolerance=500
\tolerance=1000
\brokenpenalty=5000
\frenchspacing
\catcode`\;=\active
\def;{\relax\ifhmode\ifdim\lastskip>\z@
\unskip\fi\kern.2em\fi\string;}
\catcode`\:=\active
\def:{\relax\ifhmode\ifdim\lastskip>\z@\unskip\fi
\penalty\@M\ \fi\string:}
\catcode`\!=\active
\def!{\relax\ifhmode\ifdim\lastskip>\z@
\unskip\fi\kern.2em\fi\string!}
\catcode`\?=\active
\def?{\relax\ifhmode\ifdim\lastskip>\z@
\unskip\fi\kern.2em\fi\string?}
\def\^#1{\if#1i{\accent"5E\i}\else{\accent"5E #1}\fi}
\def\"#1{\if#1i{\accent"7F\i}\else{\accent"7F #1}\fi}
\count10=0  
\count11=-1
\count12=0 
\count13=1 
\count14=1
 
\def\parno{\count11}  
\def\propno {\count12} 
\def\spropno{\count10}

\def\parag#1{\advance\parno by 1\vskip 5mm plus15mm minus 1.5mm\penalty-50
\noindent {\bf\S \the\parno-{\title {#1}}}  
\count12=0\nobreak\vs{5} } 

\long\def\sec#1#2{\count10=0\advance \propno by 1{{ \bf#1
\bf\the\parno.\bf\the\propno.
 }{\it#2}}\vskip 3mm}

\long\def\subsec#1#2{\global\advance\spropno by 1{{
\bf#1 \bf\the\parno.\bf\the\propno.\the\spropno. } {\it #2}}\vskip 3mm}

\newtoks\Title       \Title={\hfil}
\def\title#1{\Title={#1}{\leftskip=2cm plus 1fil%
     \rightskip=2cm plus 1fil%
     \parfillskip=0pt\noindent{
     \bf\the\Title\par}}}
\def\P#1#2{{\bb P}_{#1}^{#2}}
\def\A#1#2{{\bb A}_{#1}^{#2}}

\def\dem{{\it D\'emonstration. }}
 \def\ind{\hskip 1cm\relax}
\def\hs#1{\hskip #1mm}
\def\vs#1{\vskip #1mm}

\def\isom{\mathop{\fh^{\dis\hs{-5}\sim}}}

\def\liminv{\mathop{\oalign{lim\cr
\hidewidth$\longleftarrow$\hidewidth\cr}}}


\def\ie{{\it i.e.,\/}\ }

\def\up#1{\raise 1ex\hbox{\sevenrm#1}}

\def\system#1{\left\{\null\,\vcenter{\openup1\jot\m@th
\ialign{\strut\hfil$##$&$##$\hfil&&\enspace$##$\enspace&
\hfil$##$&$##$\hfil\crcr#1\crcr}}\right.}
\def\diagram#1{\def\normalbaselines{\baselineskip=20pt
\lineskip=2pt\lineskiplimit=3pt}   \matrix{#1}}


\def\imm#1{\smash{\mathop{{\hookrightarrow}}
\limits^{\scriptstyle#1}}}
\def\hfl#1#2{\smash{\mathop{\hbox to 6mm{\rightarrowfill}}
\limits^{\scriptstyle#1}_{\scriptstyle#2}}}
\def\bigver{\Big\Vert}
\def\vfl#1#2{\llap{$\scriptstyle#1$}\left\downarrow
\vbox to 4mm{}\right.\rlap{$\scriptstyle#2$}}
\def\fh{\hbox to 6mm {\rightarrowfill}\hskip 1mm}

\def\coh#1#2#3{H^{#1}({#2},{#3})} 
\def\cohr#1#2{H^{#1}({#2})}

\def\Torr#1#2#3#4{{\rm
Tor\/}^{#1}_{#2}({#3},{#4})}                                               
                                       
\def\dcoh#1#2#3{h^{#1}({#2},{#3})}
\def\dcohr#1#2{h^{#1}({#2})}

\def\fib{\mathop{\times}\limits}
\def\fibp#1#2#3{{#1}\fib_{#2} {#3}}
\def\ten{\mathop{\otimes}\limits}
\def\tens#1#2#3{{#1}\ten_{#2} {#3}}

\font\tenbb=msbm10 at 10pt       
\font\sevenbb=msbm10 at 7pt
\font\fivebb=msbm10 at 5pt        
\newfam\bbfam
\textfont\bbfam=\tenbb
\scriptfont\bbfam=\sevenbb
\scriptscriptfont\bbfam=\fivebb
\def\bb{\fam\bbfam\tenbb}

\def\hexnumber@#1{\ifcase#1 0\or1\or2\or3\or4\or5
\or6\or7\or8\or9\or 
	A\or B\or C\or D\or E\or F\fi}  
\edef\bb@{\hexnumber@\bbfam}  
\mathchardef\subsetneq="3\bb@28

\font\tenee=eufm10 at 10pt        
\font\sevenee=eufm10 at 7pt
\font\fiveee=eufm10 at 5pt        
\newfam\eefam
\textfont\eefam=\tenee
\scriptfont\eefam=\sevenee
\scriptscriptfont\eefam=\fiveee
\def\ee{\fam\eefam\tenee}
                       

\def\AA{{\cal A}}       \def\KK{{\cal K}}
      \def\LL{{\cal L}}
      \def\MM{{\cal M}}
      \def\NN{{\cal N}}
\def\EE{{\cal E}}     \def\OO{{\cal O}}
\def\FF{{\cal F}}     \def\PP{{\cal P}}
     \def\QQ{{\cal Q}}
\def\HH{{\cal H}}      
     
\def\JJ{{\cal J}}      \def\TT{{\cal T}}
\def\UU{{\cal U}}      \def\VV{{\cal V}}
\def\WW{{\cal W}}      
\def\ZZ{{\cal Z}}

\def\mod#1{\ \ ({\rm mod.\ }#1)}
\let\hook=\hookrightarrow
\let\w=\widetilde 
\def\sc{\scriptstyle}
\def\scsc{\scriptscriptstyle}
\def\dis{\displaystyle}
\def\lt{\w L }
\def\ct{\w C }
\def\nstct{\MM\hs{-1}\mid_{\ct}}
\def\ncp{\NN_{C/{{\bf P}^3}}}
\def\nctp{\MM_{\ct}}
\def\nctst{\NN_{{\ct}/{\st_0}}}
\def\ost{{{\cal O}_{\st_0}}(4H-{\w L_0}) }
\def\cq{\oplus_{i=1}^4  k(R_i)}
\def\Hdg{H_{d,g}}

\def\fitp{{\varphi^*}\TT_{\P{}{3}}}
\def\tst{\TT_{\st_0}}
\def\nstpt{\NN_{\st_0/{\w\P{}{}}}}

\def\pt{{\w \P{}{}}}
\def\st{\w S}

\def\pitp{{\pi^*}\TT_{\P{}{3}}}
\def\tptst{\TT_{\w {\P{}{}}}\hskip -1mm \mid_{\st_0}}

\def\tpt{\TT_{\w {\P{}{}}}}

\def\vs#1{\vskip #1mm}
\def\hs#1{\hskip #1mm}
\def\af#1{{\bb A}^3_{#1}}
 
\def\wpol{{\w{\bf Pol}}} 
\def\pol{{\bf Pol}}
\def\spec{{\rm Spec\ }}
\def\wp{\w {\P{}{3}}}

\def\diff{\Omega_{X\!/\!Y}}
\def\pnd#1{\P{D_{#1}}{3}} 
\def\wpnd#1{\fibp{\w {\P{}{3}}}{}{D_{#1}}}
\def\eps#1#2{\varepsilon_{ #1}^{#2}}

\def\deg{\rm deg \/}

\def\ref#1#2#3#4{{\bf[#1]}{ #2, }{\it #3, }{#4.}\par }
\def\cqfd{\unskip\kern 6pt\penalty 500
\raise -2pt\hbox{\vrule\vbox to5pt{\hrule width 4pt      
\vfill\hrule}\vrule}\vs{1}}
\def\decale#1{\par\noindent\hskip 3em\llap{#1\enspace}\ignorespaces}

\catcode`\@=12


\title{Exemples de composantes irr\'eductibles et non r\'eduites du
sch\'ema de Hilbert des courbes
lisses et connexes de $\P{}{3}$ (II)}\vs{2}
\centerline {par Sa{\"i}d AZZIZ}\vs{6}

\parag{Introduction}
\ind Soit $\Hdg$ le sch\'ema de Hilbert des courbes lisses et connexes de
degr\'e $d$ et genre $g$
dans l'espace projectif  $\P{k}{3}$ ( qu'on notera simplement $\P{}{3}$
tant que le
contexte reste clair) avec $k$ un corps alg\'ebriquement clos de
caract\'eristique nulle. Le
r\'esultat principal de cet article est le\vs{4}  
{\bf Th\'eor\`eme. }{\it Soient $d$ et  $g$ deux
entiers satisfaisant:\smallskip
$(A)$\hs{30}\ \ $G(d,8)<g\leq {\dis {73\over 2}}(d-74) $\ \  si\ \ $146\geq d\geq
95$\smallskip
$(B)$\hs{30} \ \  $G(d,8)<g\leq {\dis{(d-1)^2\over 8}}$\ \  \ si\ \
$d\geq147$\smallskip  Alors il
existe une composante irr\'eductible non r\'eduite du sch\'ema de Hilbert
$\Hdg$ dont l'\'el\'ement
g\'en\'eral est une courbe trac\'ee sur une surface quartique \`a droite
double.}\vs{4}

Ici, pour tout entier $s$, $G(d,s)$ d\'esigne le maximum des genres des
courbes lisses et connexes de
degr\'e $d$ qui ne sont pas trac\'ees sur une surface de degr\'e $s-1$;
d'apr\`es ([GP1],
th\'eor\`eme 3.1), si
$d>s(s-1)$, alors $$G(d,s)= { \dis 1\over
\dis{2s}}d^2+{ \dis s-4\over
\dis{2}}d+1-{\dis (s-1)\nu(s-\nu)\over
\dis{2s}}$$
 o\`u $\nu$ est l'unique entier tel que $0\leq \nu <s$ et $\nu\equiv d\
({\rm mod}\ s)$.\vs{1}
\ind Le premier exemple de composantes irr\'eductibles non r\'eduites fut
donn\'e par Mumford
([M1]), suivirent les exemples de Gruson et Peskine ([GP2]), Kleppe ([K1],
[K2]),
Ellia et Fiorentini ([EF]), Ellia ([E]), Hartshorne ([H1]); dans  ces
exemples les courbes
consid\'er\'ees sont situ\'ees sur des surfaces cubiques lisses. Dans
([DP])
Dolcetti et Pareschi donnent un exemple o\`u les courbes sont sur une
surface quartique \`a conique
double. Dans la majorit\'e de ces travaux, les $\Hdg$ consid\'er\'es sont
tels que $g>G(d,5)$;
notre r\'esultat en est une am\'elioration puisque $G(d,5)>G(d,8)$ pour
$d\geq 95$.    \par
\ind Le th\'eor\`eme ci-dessus ( qui est prouv\'e dans la troisi\`eme
section de cet article
) constitue en fait une application de r\'esultats plus g\'en\'eraux,
d\'emontr\'es dans la
deuxi\`eme partie, qui concernent les d\'eformations des courbes de
$\P{}{3}$
situ\'ees sur des surfaces
\`a singularit\'es ordinaires de lieu double une courbe lisse. Notre
motivation est le fait que les
surfaces quartiques
\`a droite double que nous consid\'erons sont en particulier des surfaces
\`a singularit\'es
ordinaires. \par
\ind Le plus difficile dans la preuve de ce th\'eor\`eme et de montrer que
si $C\in\Hdg$ est une
courbe trac\'ee sur une quartique \`a droite double, alors toute
g\'en\'erisation de
$C$ dans $\Hdg$ est aussi trac\'ee sur une surface quartique \`a droite
double. On utilisera pour 
cela la proposition (2.11) qui tra\^ite de la m\^eme question mais plus
g\'en\'eralement pour les
courbes trac\'ees sur des surfaces \`a singularit\'es ordinaires.\par
\ind La premi\`ere partie de l'article rappelle les  d\'efinitions des
sch\'emas polaires
et des pinceaux de Lefschetz; les r\'esultats qui y sont \'etablis
constituent les pr\'eparatifs
n\'ecessaires pour \'etablir ceux de la seconde partie qui, elle, est la
pierre angulaire du
pr\'esent article.\vs{5}
\ind Cet article  g\'en\'eralise les r\'esultats d'un travail de th\`ese
(cf. [Az]) effectu\'e 
sous la direction de N. Mestrano, qu'elle trouve ici mes sinc\`eres
remerciements. Je
remercie aussi Carlos Simpson pour les nombreuses heures de discussion
qu'il
m'a accord\'e et o\`u j'ai beaucoup appris.
\parag{Sch\'emas polaires } Dans cette section nous traitons des {\it
sch\'emas polaires}, introduits
par B. Teissier (cf. [T1], [T2]) et D.S. Rim (cf. [R]). Nous en donnons une
construction inspir\'ee
de ces travaux, ensuite nous \'etablissons le lien existant entre des
sch\'emas polaires
particuliers et la notion de pinceaux de Lefschetz (cf. 1.7). Nous
terminerons par la proposition
(1.10) qui sera l'outil essentiel utilis\'e dans la preuve de la
proposition-cl\'e de la prochaine
section (cf. 2.5).\vs{3} 
\ind Soient $X$ et $Y$ deux $k$-sch\'emas  irr\'eductibles de type fini,
soit $f\colon X\to Y$ un
morphisme lisse de dimension relative $m$, soit
$\LL$ un fibr\'e en droites sur $X$. Soit $\diff(\!=\!\Omega^1_{X\!/\! Y})$
le fibr\'e des
diff\'erentielles relatives de $X$ sur $Y$ muni du morphisme
$d\colon\OO_X\fh\diff$.\par
\ind Soit $P$ un point de $X$, on dit qu'un voisinage affine $U$ de $P$ est
un {\it voisinage de
coordonn\'ees locales relatives}, s'il existe
$x_1,\dots ,x_m\in \OO_X(U)$ tels que $x_i(P)\hskip -0.5mm=\hskip -0.5mm
0$, et
$\diff\hskip -1.5mm\mid_U$$=\oplus_{i=1}^m\OO_Udx_i$ (de tels voisinages existent car $f$ est
lisse).\par   
\ind Soit $(U_\alpha )_{\alpha \in \AA}$ un recouvrement ouvert de $X$
form\'e par des voisinages
affines de coordonn\'ees locales relatives $(x_1^\alpha,\dots
,x_m^\alpha)_\alpha$ tel qu'il existe
une trivialisation $\theta\colon\LL\fh\OO_X$ de $\LL$ par rapport \`a
$(U_\alpha )_\alpha$, \ie, une
famille d'isomorphismes
$(\theta_{\alpha}\colon\LL\mid_{\dis {U_\alpha}}\!\isom\
\OO_{U_\alpha})_{\alpha}$. Soit $\sigma$ une section globale non nulle de
$\LL$. Pour tout
$\alpha\!\in\!\AA$, notons $
\sigma_{\alpha}=\theta_{\alpha}(\sigma\mid_{U_{\alpha}})$ et consid\'erons
l'id\'eal\smallskip
$\hs{35}I_\alpha=\bigl(\dis\sigma_{\alpha },{{\partial \sigma_{\alpha
}}\over{\partial x_1^{\alpha
}} },\dots ,{{\partial \sigma_{\alpha }}\over{\partial x_m^{\alpha }}
}\bigr)\subset\OO_{U_\alpha}(U_\alpha)$\smallskip Soit $\JJ_\alpha=\w
I_\alpha$ le faisceau
d'id\'eaux de $\OO_{U_\alpha}$ associ\'e
\`a $I_\alpha$, alors les faisceaux $\JJ_\alpha $ se recollent le long des
$U_\alpha$ en un faisceau d'id\'eaux de $\OO_X$, que l'on notera $\JJ$. En
effet pour tous $\alpha$
et $\beta$ dans $\AA$, si $U_{\alpha\beta}$ d\'esigne
$U_{\alpha}\cap U_{\beta}$, on a
$\sigma_{\alpha\beta}=t_{\alpha\beta}\sigma_{\beta\alpha}$,  o\`u
$\sigma_{\beta\alpha
}\!=\!\sigma_{\beta}\!\mid_{U_{\alpha\beta}}$,
$\sigma_{\alpha\beta }\!=\!\sigma_{\alpha}\!\mid_{U_{\alpha\beta}}$ et
$t_{\alpha\beta}\in\coh{0}{U_{\alpha\beta}}{\OO_X^*}$, ainsi \smallskip
\hs{50}$\dis{{\partial \sigma_{\alpha\beta }}\over{\partial
x_i^{\alpha\beta }}
}=t_{\alpha\beta}{{\partial \sigma_{\beta\alpha }}\over{\partial
x_i^{\beta\alpha }}
}+\sigma_{\beta\alpha }{{\partial t_{\alpha\beta}}\over{\partial
x_i^{\beta\alpha }} }$\smallskip
 ce qui implique
$\JJ_\alpha\!\!\mid_{U_{\alpha\beta}}=\JJ_\beta\!\!\mid_{U_{\alpha\beta}}$
comme sous-faisceaux de
$\OO_{U_{\alpha\beta}}$. Un calcul analogue permet de voir que le
sous-sch\'ema correspondant au
faisceau  d'id\'eaux $\JJ$ ne d\'epend pas du recouvrement
$\{U_\alpha\}_{\alpha}$ ni de la
trivialisation
$\{\theta_{\alpha}\}_{\alpha}$ choisis. Ceci donne alors un sens \`a la
d\'efinition suivante:\vs{2}
\sec{D\'efinition}{ Soient $X$, $Y$, $f$, $\LL$, $\sigma$ et $\JJ$ comme
ci-dessus. On appelle
sch\'ema polaire de $\sigma$ par rapport au morphisme
$f$, le sous-sch\'ema ferm\'e de $X$ correspondant \`a $\JJ$ qu'on notera
\break  $ \pol(f\colon
X\to Y,\sigma)$ ou simplement ${\bf Pol}(f,\sigma)$ quand le contexte est
clair.}
\ind D'apr\`es la construction faite ci-dessus on peut voir sans
difficult\'e la\vs{1}
\sec{Proposition}{ Soit $X$, $Y$, $f$, $\LL$ et
$\sigma$ comme ci-dessus. Alors \par {\parindent=1cm 
\item{\rm (i)} $\pol(f,\sigma)$ est \'egal au lieu singulier du morphisme
$f\!\mid_{\Sigma}$ o\`u $\Sigma$ est le sch\'ema des z\'eros de $\sigma$.
\item {\rm (ii)} Pour tout changement de base $h\colon Y'\to Y$ on a ${\bf
Pol}(\fibp{X}{Y}{Y'}\to
Y',q^*\sigma)=
\fibp {\pol(f,\sigma)}{Y}{Y'}$ o\`u $q$ d\'esigne la projection
$\fibp{X}{Y}{Y'}\to X$.\par}}
\ind Dans ce qui suit, nous allons d\'ecrire les sch\'emas polaires des
surfaces de $\P{}{3}$
relativement aux morphismes induits par l'\'eclatement de
$\P{}{3}$ le long d'une droite g\'en\'erale, en \'etablissant le lien qui
existe entre ces
sch\'emas polaires et la notion des pinceaux de Lefschetz.\vs{2}
\sec{D\'efinitions}{ Soit ${\P{}{3}}^{\vee}$ l'espace des plans de
$\P{}{3}$. Un point $\Pi$ de
${\bb G}(1,{\P{}{3}}^{\vee}) $, la grassmannienne des droites de
${\P{}{3}}^{\vee}$, s'appelle un pinceau de plans dans $\P{}{3}$.
Consid\'erant $\Pi$ comme une
droite dans
${\P{}{3}}^{\vee}$, on \'ecrira par abus de notation
$\Pi\!=\!{\{H_t\}}_{t\in\P{}{1}}$. On appelle axe du pinceau  $\Pi$
l'intersection de deux plans
distincts quelconques de $\Pi$.}  
\sec{D\'efinition}{ Soit $S\subset\P{}{3}$ une surface lisse, un pinceau 
 $\Pi={\{H_t\}}_{t\in\P{}{1}}$ s'appelle un pinceau de Lefschetz de $S$
s'il v\'erifie les axiomes
suivants:\par  {\rm L1)} L'axe $\Delta$ de $\Pi$ coupe transversalement
$S$.\par    {\rm L2)} Il
existe un ouvert dense $\!U$ de $\P{}{1}$ tel que pour tout $t\!\in\! U$,
$H_t$ coupe
transversalement $S$, et\par  {\rm L3)} pour tout
$t_0\!\in\!\P{}{1}\setminus U $,
$H_{t_0}$ coupe transversalement $S$ sauf en un point $x_0$, qui est un
point singulier quadratique
ordinaire de $H_{t_0}\cap S$, i.e., un point de multiplicit\'e 2 avec des
directions tangentes
distinctes. (Le point $x_0$ est appel\'e point singulier du pinceau $\Pi$,
on notera ${\rm Sing\/}
\Pi$ l'ensemble de tels points.)}
\ind La notion des pinceaux de Lefschetz a \'et\'e trait\'ee par N.M.Katz
(cf. [K1], [K2]) dans un
cadre plus g\'en\'eral qui est celui des sou-sch\'emas lisses et
irr\'eductibles de $\P{}{n}$ pour
$n$ quelconque. Les r\'esultats de ces travaux nous concernant sont
r\'eunis dans le \smallskip
\sec{Th\'eor\`eme}{Soit $S$ une surface lisse de $\P{}{3}$. Alors\par
 \hs{9}{\rm (i)} les pinceaux de Lefschetz de $S$ existent {\rm([Ka1],
2.5)}, \par
\hs{9}{\rm (ii)} ces pinceaux forment un ouvert dense de
${\bb G}(1,{\P{}{3}}^{\vee})$ {\rm([Ka1], 3.2.1)}.\par
\hs{9}{\rm (iii)} Si $\Pi$ est un pinceau de Lefschetz de $S$, alors le
cardinal de Sing$\Pi$ (qui
est fini par l'axiome {\rm L3}) ne d\'epend pas de $\Pi$, plus
pr\'ecis\'ement, c'est le degr\'e de
la surface duale de $S$, i.e., $\#{\rm Sing\ }\Pi=s(s-1)^2$  avec $s=\deg\
S$ {\rm([Ka2],
3.2.10)}.}\smallskip
\sec{Remarque}{\rm Si $\Pi$ est un pinceau de Lefschetz d'une surface
$S\subset\P{}{3}$, d'axe $\Delta$, alors les axiomes (L1) et (L3)
entra\^{\i}nent le fait que
Sing$\Pi$ est contenu dans $S\setminus\Delta$. On montre, dans le
r\'esultat qui suit, que Sing
$\Pi$ est lisse et correspond \`a un sch\'ema polaire particulier.}
\sec{Th\'eor\`eme}{ Soit $S$ une surface lisse de $\P{}{3}$ de degr\'e $s$,
d\'efinie par un
polyn\^ome homog\`ene $F$. Soit $\Pi$ un pinceau de Lefschetz de
$S$ d'axe $\Delta$. Soit $\eta\colon\wp\to\P{}{3}$ l'\'eclatement de
$\P{}{3}$ en
$\Delta$, et soit $q\colon\wp\to\P{}{1}$ la projection induite sur
$\P{}{1}$. Notons
$\w F=\eta^*F$, alors \par
\hskip 1cm{\rm (i)} $\pol(q,\w F)=\eta^{-1}({\rm Sing} \Pi)$.\par
\hskip 1cm{\rm (ii)} $\pol(q,\w F)$ est fini de cardinal \'egal \`a
$s(s-1)^2$.\par
\hskip 1cm {\rm (iii)} Tous les points de $\pol(q,\w F)$ sont simples.}
\dem Tout d'abord, comme $\eta$ est un isomorphisme en dehors de $\Delta$,
on voit que l'assertion
(ii) est une cons\'equence de (i) et de la remarque (1.6).\par
\ind Soit $\{X_0,X_1,X_2,X_3\}$ un syst\`eme de coordonn\'ees homog\`enes
de
$\P{}{3}$ tel que la droite $\Delta$ soit donn\'ee par ${X_0}=X_1=0$.
Notons
$\af{i}=\P{}{3}\setminus\{X_i=0\}$, et $\w {\af{i}}=\eta^{-1}({\af{i}})$.
On a par d\'efinition
\smallskip
$\pol(q,\w F)={\bigcup^3_{i=0}}Pol_i$ avec 
$Pol_i=\pol(q\mid_{\w\af{i}},\w F_i\!=\!\w F\mid_{\w {\af{i}}})$, ainsi
pour prouver les assertions
(1) et (3), il nous suffit de d\'eterminer les $Pol_i$ pour tout $i$.\par
$ {\underline{1^e\ Cas: i\in\{2,3\}}}$. Supposons que $i=2$ (le cas $i=3$
se d\'emontre de la m\^eme
mani\`ere). Soient
$x\!=\! X_0/ X_2,\ y\!=\!X_1/X_2$, et $z=X_3/X_2$ les coordonn\'ees affines
de $\af{2}$, alors la
droite $\Delta$ est d\'efinie par
$x\!=\! y\!=0$ et l'ouvert $\w \af{2}$ est r\'eunion des deux ouverts
affines $
{\af{(x,w,z)}}\!=\!\spec k[x,w,z]	$ et ${\af{	(y,v,z)}}=\spec k[v, y,z]$,
o\`u
$w=y/x$ et $v=x/y$. Le morphisme $\eta\!\mid_{\af{(x,w,z)}}$ (resp.
$\eta\!\mid_{\af{(y,v,z)}}$) s'\'ecrit $(x,w,z)\mapsto\ (x,xw,z)$ ( resp.
$(y,v,z)\mapsto
(vy,y,z))$, et le morphisme $q\!\mid_{\af{(x,w,z)}}$ (resp.
$q\!\mid_{\af{(y,v,z)}}$) s'\'ecrit $(x,w,z)\mapsto w$ (resp.
$(y,v,z)\mapsto v)$. Soit
$F_2(x,y,z)=F(x,y,1,z)$ l'\'equation de $S$ dans $\af{2}$, alors on a 
$$Pol_2 =V(\w {F_2}, {\partial  \w {F_2}\over{\partial  x }}, {\partial
\w {F_2}\over{\partial  z}})	 \cup V({\w{F_2}}',{\partial {\w
{F_2}}'\over{\partial  y }},
{\partial{\w {F_2}}'\over{\partial z}})$$	 o\`u $\w
{F_2}(x,w,z)={F_2}(x,xw,z)$, et
$\w{F_2}'(y,v,z)={F_2}(yv,y,z)$.
 Soit $M_0=(x_0, w_0, z_0)$ un point de $V(\w {F_2},{\dis{\partial
\w {F_2}\over{\partial  x }},{\partial \w {F_2}\over{\partial  z }}})	$ ayant pour image
$N_0=(x_0, w_0x_0,z_0)$ et consid\'erons le plan $H_{(w_0,1)}$ d'\'equation
\smallskip 
$w_0x\ -\ y=0$, on a alors $H_{(w_0,1)}\in \Pi$ et $N_0$ est le n{\oe}ud de
la courbe $
S_{(w_0,1)}=H_{(w_0,1)}\cap S$, i.e., $N_0\in$ Sing$\Pi$. En effet la
matrice jacobienne de
$S_{(w_0,1)}$ a pour mineurs au point $N_0$ les quantit\'es
$$\alpha={\dis{\partial F_2\over{\partial  x }}\mid_{N_0}+w_0{\partial 
F_2\over{\partial
y}}\mid_{N_0}}\ ,\ \beta=w_0{\partial \dis F_2\over{\partial
\dis z }}\mid_{N_0}\ ,\ \gamma={\partial \dis F_2\over{\partial \dis z
}}\mid_{N_0}$$
\vs{-3}qui sont nuls car $\alpha={\dis{\partial \w {F_2}\over{\partial x
}}\mid_{M_0}}=0$,
$\beta=w_0{\dis{\partial \w {F_2}\over{\partial z }}\mid_{M_0}}=0$ et
$\gamma={\dis{\partial
\w {F_2}\over{\partial z }}\mid_{M_0}}=0$. Ce raisonnement \'etant valable
si l'on remplace $\w
{F_2}$ par ${\w {F_2}}'$ on a donc $Pol_2\!\subset\!\eta^{-1}($Sing$\Pi)$.
Inversement soit
$N_0\!=\!(x_0,y_0, z_0)$ un point de Sing$\Pi$, alors $(x_0,y_0)\not
=(0,0)$ car
$($Sing$\Pi)\cap\Delta=\emptyset$; soit par exemple $x_0\not =0$, on a donc
$N_0=\eta (x_0,y_0/x_0,z_0)$ et $N_0$ ( d'apr\`es le calcul pr\'ec\'edent)
est dans
 $V(\w {F_2}, {\dis{\partial \w {F_2}\over{\partial x }}, {\partial\w
{F_2}\over{\partial z }}})$,
i.e., dans $Pol_2$. Ceci d\'emontre l'\'egalit\'e de (i). Pour prouver
(ii), il faut montrer que le
d\'eterminant de  de la matrice jacobienne $J_2=$Jac$_{(x,w,z)}(\w
{F_2},{\dis{\partial
\w {F_2}\over{\partial  x }},{\partial \w {F_2}\over{\partial  z }}})$
(resp.
$J_2'=$Jac$_{(y,v,z)}({\w {F_2}}',{\dis{\partial {\w {F_2}}'\over{\partial 
y }},{\partial {\w
{F_2}}'\over{\partial  z }}})$) est non nul en tout point $M$ de
$V(\w{F_2},{\dis{\partial
\w{F_2}\over{\partial x}},{\partial \w{F_2}\over{\partial z }}})$ (resp.
$V({\w {F_2}}',{\dis{\partial {\w {F_2}}'\over{\partial y}},{\partial {\w
{F_2}}'\over{\partial z
}}})$). Un calcul simple montre que
det$J_2\mid_M={\dis{\partial\w{F_2}\over{\partial 
w}}\biggl({\partial^2 
\w{F_2}\over{\partial  x^2 }}{\partial^2 
\w{F_2}\over{\partial  z^2 }}-({\partial^2  \w{F_2}\over{\partial 
x\partial
z}})^2\biggr)}\mid_{M}$, et comme ${\dis\partial\w{F_2}\over{\dis\partial 
w}}\mid_{M}\not=0$ (car
sinon $\w S=\eta^{-1}S$ serait singuli\`ere), alors il suffit de voir que
$A=\biggl({\dis\partial^2  \w{F_2}\over{\dis\partial  x^2 }}\
{\dis\partial^2  \w
F_2\over{\dis\partial  z^2 }}\ -\ ({\dis\partial^2 
\w{F_2}\over{\dis\partial  x\dis\partial
z}})^2\biggr)\mid_{M}\not =0$.\par  Pour cela on consid\`ere un point
$M_0=(x_0,w_0,z_0)$ de
$V(\w{F_2},{\dis{\partial
\w{F_2}\over{\partial x}},{\partial \w{F_2}\over{\partial z }}})$ d'image
$N_0$ par $\eta$. On a
$x_0\not=0$ puisque $x_0{\dis\partial F_2\over{\dis\partial
y}}\mid_{N_0}={\dis\partial 
\w{F_2}\over{\dis\partial  w}}\mid_{M_0}\not=0$, donc
$S_{(w_0,1)}$ n'est pas contenue dans le plan $x=0$ et, si $S'_{(w_0,1)}$
d\'esigne la projet\'ee de
$S_{(w_0,1)}$ sur le plan $y=0$, on a $S_{(w_0,1)}\simeq S'_{(w_0,1)}$ et
$N'_0=(x_0,w_0)$ est le
n{\oe}ud de $S'_{(w_0,1)}$. Soit $h(x,z)=F_2(x,w_0x,z)$ l'\'equation de la
courbe $X'_{(w_0,1)}$, le
fait que le point $N'_0$ soit un n{\oe}ud \'equivaut \`a
$B=({\dis\partial^2 h\over{\dis\partial x^2
}}{\dis\partial^2 h\over{\dis\partial  z^2 }}\ -\ ({\dis\partial^2 
h\over{\dis\partial 
x\dis\partial z}})^2)\mid_{N_0'}\not =0$, et il est clair que
$A\!=\!B$.\par 
  ${\underline {2^e\ Cas: \ i\in\{0,1\}}}$ Prenons $i=0$. Comme ${\w
{\af{0}}}\subset
(\wp\setminus\eta^{-1}(\Delta))\simeq
\P{}{3}\!\setminus\! \Delta$, alors le morphisme $q\mid_{\w
{\af{0}}}$\break n'est autre que la
restriction de la projection de centre $\Delta$, i.e., la fl\`eche
$(X_0,X_1,X_2,X_3)\!\mapsto{(X_0,X_1)}$.\par Soient $x=X_1/X_0,\
y=X_2/X_0,\ z=X_3/X_0$ les
coordonn\'ees affines dans
$\af{0}$. On a $Pol_0=V(F_0,{\dis {\partial F_0\over{\partial y }},
{\partial  F_0\over{\partial z
}}})$ avec $F_0(x,y,z)= F(1,x,y,z)$, et puisque l'\'equation de tout
\'el\'ement de $\Pi$ est de la
forme
$aX_0+bX_1=0$ o\`u $(a,b)\in \P{}{1}$, alors on voit qu'un point $M$ est dans
$Pol_0$ si et seulement si le plan tangent $T_MS$ est dans $\Pi$ (ici $M$
est identifi\'e \`a un
point de $\P{}{3}$), d'o\`u l'assertion (i).\par Pour voir (ii) il suffit
de montrer que le
d\'eterminant de la matrice $J_0=$Jac$_{(x,y,z)}( F_0,{\dis{\partial
F_0\over{\partial  x
}},{\partial F_0\over{\partial  z }}})$ est non nul sur $Pol_0$. Soit $M\in
Pol_0$, on a
det$J_0\mid_M={\dis{\partial F_0\over{\partial x}}({\partial^2 
F_0\over{\partial y^2 }}{\partial^2 
F_0\over{\partial  z^2 }} - ({\partial^2  F_0\over{\partial  y\partial
z}})^2)}\mid_{M}$
 et comme ${\dis{\partial F_0\over{\dis\partial x}}}\mid_{M}\not =0$ (sinon
$S$ serait
singuli\`ere), alors il suffit de voir que
$A=({\dis{\partial^2  F_0\over{\partial y^2 }} {\partial^2 
F_0\over{\partial  z^2 }}-({\partial^2 
F_0\over{\partial  y\partial z}})^2)}\mid_{M}\not =0$. En raisonnant comme
dans le premier cas, on
voit que $A\not =0$ \'equivaut au fait que le point $M$ est de
multiplicit\'e 2 dans la courbe
$S\cap T_MS$, mais ceci est assur\'e par l'axiome (L3) (cf. 1.4).\cqfd
\vs{1}
\ind La preuve de la proposition ci-dessus montre que si l'on sait
d\'efinir les pinceaux de
Lefschetz d'une surface singuli\`ere, ( c'est le cas pour les surfaces
ayant un nombre fini de
singularit\'es, et \`a ma connaissance c'est le seul), on pourrait alors
g\'en\'eraliser
ladite proposition \`a une surface $S$ singuli\`ere et irr\'eductible; dans
ce cas le sch\'ema des
p\^oles serait isomorphe
\`a Sing $S\ \cup\ $Sing$\ \Pi$ o\`u $\Pi$ est un pinceau de Lefschetz pour
$S$. Faute de pouvoir
d\'efinir les pinceaux de Lefschetz d'une surface singuli\`ere, le
th\`eor\`eme (1.7) se laisse
g\'en\'eraliser de la fa{\c c}on suivante:\vs{2}
\sec{Proposition}{Soit $S=V(F)$ une surface singuli\`ere et irr\'eductible
de $\P{}{3}$. Soit
$\Delta$ une droite disjointe de Sing $S$ coupant $S$ transversallement. On
d\'efinit deux
morphismes $\eta$ et $q$ comme dans (1.7). Alors
$\eta^{-1}($Sing $S)$ est un sous-sch\'ema ferm\'e de $\pol(q,\eta^*F)$ et
$\dim\pol(q,\eta^*F)\leq1$.}
\dem On va utiliser les notations de la preuve de (1.7). La question
\'etant locale, on se place par
exemple dans $\w\A{2}{3}$. (Le cas de $\w\A{3}{3}$ est analogue, celui de
$\w\A{0}{3}$ et
$\w\A{1}{3}$ est trivial.) Soit
$M_0=(x_0,w_0,z_0)\in\eta^{-1}($Sing $S)$, alors\par 
\hs{20}${\dis{\partial \w {F_2}\over{\partial x
}}\mid_{M_0}}={\dis{\partial F_2\over{\partial x
}}\mid_{\eta(M_0)}+w_0{\partial  F_2\over{\partial y}}
\mid_{\eta(M_0)}}=0\ {\rm et\ }{\dis{\partial
\w {F_2}\over{\partial z }}\mid_{M_0}}={\dis{\partial F_2\over{\partial z
}}}
\mid_{\eta(M_0)}=0,$\smallskip autrement dit $M_0\!\in\! Pol_2$. Pour
montrer la deuxi\`eme partie
de l'assertion, on sait que $\dim\pol(q,\eta^*F)$ vaut au plus $2$ car
$\pol(q,\eta^*F)$ est contenu
dans
$\st$ et il est facile de voir que les points de
$\eta^{-1}[S\setminus($Sing $S\cup (\Delta\cap
S))]$ \'evitent $\pol(q,\eta^*F)$, i.e.,
$\pol(q,\eta^*F)\not=\st$.\cqfd\vs{2}
\ind Nous arrivons maintenant \`a la proposition-synth\`ese de cette
section. On utilisera les
\vs{1}\sec{Notations}{\rm\ Pour tout entier $n\geq0$, on note $D_n=\spec
k[t]/_{\dis (t^{n+1})}$,
o\`u
$t$ est une ind\'etermin\'ee sur $k$. On note $\A{}{1}=\spec k[t]$. Par
$\{1\}$ (resp. $\{0\}$) on
d\'esigne le point de $\A{}{1}$ correspondant \`a l'id\'eal $(t-1)$ (resp.
$(t)$), autrement dit,
$\{1\}=\spec k[t]/_{\dis (t-1)}$ et $\{0\}=D_1$.}  
\sec{Proposition}{ Soit $\Delta$ une droite g\'en\'erale de $\P{}{3}$, on
consid\`ere
$\eta_0\colon\wp\to\P{}{3}$ l'\'eclatement de $\P{}{3}$ en
$\Delta$ et $q_0\colon\wp\to\P{}{1}$ la projection induite sur $\P{}{1}$.
Soit
$s$ un entier positif, on consid\`ere
$F$ un \'el\'ement de $\cohr{0}{\OO_{\P{}{3}\times\A{}{1}}(s)}$ tel 
$G=F\mid_{\P{}{3}\times \{1\}}$ est lisse et $F_n=F\mid_{\P{}{3}\times
D_n}$ est non nul pour tout
$n$. On pose  $\w\pol=\pol(q_0\times id_{\A{}{1}},(\eta_0\times
id_{\A{}{1}})^*F)$ et
$\pol_n=\pol(q_0\times id_{D_n},(\eta_0\times id_{D_n})^*F_n)$. Alors\par
{\rm(i)} Pour tout $n\geq
0$, $\pol_0=\fibp{\pol_n}{D_n}{\spec k}$ et
$\pol_n=\fibp{\w\pol}{\A{}{1}}{D_n}=\w\pol\cap (\wp\times D_n)$.\par
{\rm(ii)}
$\fibp{\w\pol}{\A{}{1}}{\{1\}}=\wpol\cap (\wp\times\{1\})=\pol(q_0,\eta_0^*
G)$, ce dernier est un
sch\'ema fini lisse de cardinal $l(s):=s(s-1)^2$.\par {\rm(iii)} Le nombre
de points de la fibre
sp\'eciale (i.e., au-dessus de
$0$) de la projection $\wpol'\to\A{}{1}$ est au plus $l(s)$, o\`u
${\wpol}'$ d\'esigne la fermeture
dans $\wp\times\A{}{1}$ de $\wpol\setminus\pol_0$.\par
\ind On suppose de plus que {\rm Sing }$V(F_0)$ est une courbe, alors\par
{\rm(iv)} $\w\pol$ est un
sous-sch\'ema ferm\'e de $\wp\times\A{}{1}$ localement intersection
compl\`ete de dimension $1$; en
particulier $\w\pol$ est Cohen-Macaulay.\par {\rm(v)} Si $f$ est une
fonction non nulle d\'efinie
sur $\wpol$, alors (a) la dimension de toute composante irr\'eductible de
{\rm Supp}$(f)$ (i.e., le
support de $f$) est $\geq1$ et (b) tout point isol\'e de {\rm
Supp}$(f)\cap({\w\P{}{3}}\times\{0\})$
est dans
$\wpol'$; le nombre de tels points est $\leq l(s)$. 
\par} 
\dem (i) C'est clair d'apr\`es (1.2, (ii)).\par
\ind (ii) La premi\`ere partie est une cons\'equence de (1.2, (ii)). Pour
prouver la deuxi\`eme
partie, et en vue d'appliquer le tho\'er\`eme (1.7), il suffit de voir que
l'on peut consid\'erer
$\Delta$ comme l'axe d'un pinceau de Lefschetz de la surface
$V(G)$. Ceci est possible car, d'une part $\Delta$ est g\'en\'erale dans
${\bb G}(1,\P{}{3})$ et on a un isomorphisme ${\bb
G}(1,{\P{}{3}}^{\vee})\isom
\ {\bb G}(1,\P{}{3})$ d\'efini par la fl\`eche\
$\Pi\mapsto$ l'axe de $\Pi$; d'autre part, d'apr\`es (1.5, (ii)), un
\'el\'ement g\'en\'eral de ${\bb
G}(1,{\P{}{3}}^{\vee})$ peut \^etre vu comme un pinceau de Lefschetz de
$V(G)$.\par 
\ind (iii) D'apr\`es le lemme (1.10.1) ci dessous, il existe un voisinage
ouvert $W$ de
$0$ dans $\A{}{1}$ tel que le morphisme $\wpol'\cap (\wp\times W)\to W$
soit plat et fini de fibre
g\'en\'erale form\'ee de $l(s)$ points simples, ainsi la fibre sp\'eciale
est form\'ee de $m$
points, chacun de multiplicit\'e \'egale \`a $p_i\geq 1$ avec
$\dis\Sigma_1^mp_i=l(s)$, d'o\`u
$m\leq l(s)$.\par
\ind(iv)$\ $ La fibre sp\'eciale de$\wpol\to\A{}{1}$ (qui est $\pol_0$
d'apr\`es (i)) est de
dimension 1 puisqu'elle contient un sous-sch\'ema isomorphe \`a {\rm
Sing}$V(F_0)$ (cf. 1.8), et
(quitte \`a se restreindre \`a un ouvert de $\A{}{1}$ contenant $0$ et $1$)
la fibre g\'en\'erale est
de dimension 0 donc $\dim\wpol=1$. Comme $\wpol$ est d\'efini localement
par trois
\'equations dans un espace lisse de dimension 4, alors il est localement
intersection compl\`ete et
donc, en vertu de ([H2], II-8.23), il est Cohen-Macaulay.\par
\ind (v) L'assertion (a) est une cons\'equence de (1.11). Pour montrer (b),
on consid\`ere un point
$P$ isol\'e de Supp$(f)\cap (\wp\times\{0\})$. Soit $Y\subset$Supp$(f)$ une
composante irr\'eductible
contenant
$P$; d'apr\`es (a) on a $\dim Y\geq 1$, et comme $P$ est isol\'e dans
Supp$(f)\cap(\wp\times)\{0\}$
, alors $Y\not\subset\wp\times \{0\}$. Ceci implique que le point
g\'en\'eral de $Y$ est dans
$\wpol'$, or $\wpol'$ est un ferm\'e, donc $Y\subset\! \wpol'$ d'o\`u
$P\in\wpol'$.\cqfd\vs{1}
\subsec{Lemme}{ Avec les notations ci-dessus, il existe un voisinage ouvert
$W$ de
$0$ dans $\A{}{1}$ tel que le morphisme de projection
$\wpol'\cap (\wp\times W)\to W$ soit plat et fini de fibre g\'en\'erale
form\'ee de
$l(s)$ points simples.} 
\dem Par hypoth\`ese $F\mid_{\wp\times \{1\}}$ est lisse, donc il existe un
voisinage ouvert
$W'\subset\A{}{1}$ de $1$ tel que pour tout point ferm\'e $\alpha\in W'$ on
ait $F\mid_{\wp\times
\{\alpha\}}$ lisse. En utilisant (1.7) on d\'eduit que le morphisme
$\wpol\cap(\wp\times W')\to W'$
est fini de fibres form\'ees de $l(s)$ points simples, et quitte
\`a r\'eduire $W'$ on peut le supposer plat. On a $\wpol\cap(\wp\times
W')=\wpol'\cap(\wp\times W')$
(car $0\not\in W'$); on pose $W=W'\cup\{0\}$, c'est un voisinage de $0$
dans $\A{}{1}$. Quitte \`a
remplacer $\wpol'$ par $\wpol'_{\rm red}$ on peut supposer que $\wpol'$ est
r\'eduit, ainsi pour
montrer le lemme, il suffit de montrer qu'il n'y a pas de composante
irr\'eductible de
$\wpol'\cap (\wp\times W)$ dans la fibre au-dessus du point $0$ (car $W$
est un sch\'ema int\`egre
de dimension $1$). Pour ce faire, soit $\alpha\in
\wpol'\cap(\wp\times\{0\})$ un \'el\'ement de la fibre au-dessus de $0$,
alors il existe
$\beta\in\wpol\setminus(\wpol\cap
\wp\times\{0\})$ tel que $\alpha\in{\overline{\{\beta\}}}$. Le point
$\beta$ ne peut pas \^etre dans
la fibre au-dessus de $0$ (par d\'efinition), ni dans la fibre au dessus
d'un point ferm\'e $a\not
=0$ de $\A{}{1}$ car sinon $\alpha$ y serait aussi, donc $\beta$ est dans
la fibre du point
g\'en\'erique de
$\A{}{1}$, \ie $\beta\in\wpol\cap(\wp\times W')=\wpol'\cap(\wp\times W')$.
On en d\'eduit que
$\beta$ est dans une composante irr\'eductible de $\wpol'$ qui ne s'envoie
pas sur $0$; or
$\alpha\in{\overline{\{\beta\}}}$, donc $\alpha$ est dans une composante
irr\'eductible de $\wpol'$
qui ne s'envoie pas sur
$0$. Cela implique que $\alpha$ ne peut pas \^etre le point g\'en\'erique
d'une composante de
$\wpol'$ qui s'envoie sur
$0$.\cqfd\vs{2}\sec{Lemme}{Soit $X$ un sch\'ema Cohen-Macaulay. Soit
$f\in\Gamma(X,\OO_X)\setminus\{0\}$. Alors pour toute composante
irr\'eductible
$Y$ de {\rm Supp}$(f)$ on a $\dim Y\geq1$.}
\dem Soit $x$ un point de Supp$(f)$. Notons $ B=\OO_{X,x}$ et $\bf m$
l'id\'eal maximal de $B$.
Comme $B$ est un anneau Cohen-Macaulay, alors il existe $\delta\in B$ tel
que $\delta(x)=0$ et
$\delta$ est non-diviseur de z\'ero dans
$B$. On consid\`ere $I= {\rm Ann}(f_x)=\{z\in B\mid\ z.f_x=0\}$, alors $I$
est un id\'eal de $B$ tel
que $V(I)=$ Supp$(f)$ au voisinage de $x$. Supposons que
$\dim_x$ Supp $(f)=0$, alors au voisinage de $x$, Supp$(f)$ est
ensemblistement
\'egal \`a $\{x\}$; autrement dit $\sqrt I={{\bf m} }$. On a
$\delta\in{{\bf m}}$ car
$\delta(x)=0$, soit $n$ le plus petit entier tel que $\delta^n\in I$; alors
 
$0=\delta^nf_x=\delta(\delta^{n-1}f_x)$. Ceci est absurde car, d'une part
$\delta^{n-1}f_x\not =0$
par d\'efinition de $n$, et d'autre part $\delta$ est non-diviseur de
z\'ero. \cqfd

\parag{D\'eformations des surfaces \`a singularit\'es ordinaires dans
$\P{}{3}$}
\sec{}{\rm\ Dans toute cette section on d\'esigne par $S_0$ une surface
int\`egre de
$\P{}{3}$ \`a singularit\'es ordinaires, i.e., le lieu singulier de $S_0$
est une courbe double
$L_0$ lisse sauf en un nombre fini de points triples $T_1,\dots ,T_t$ et
contenant  un nombre fini
de {\it points-pinces} $P_1,\dots ,P_{p}$ (distincts des points triples)
telle que pour tout $Q\in
L_0$ on ait
$$\hat\OO_{S_0, Q}\simeq\cases{k[[t_1,t_2,t_3]]/_{\dis (t_1t_2t_3)}& si Q
$\in\{ T_1,\dots ,T_t\}$;\cr 
 k[[t_1,t_2,t_3]]/_{\dis (t_2^2-t_1^2t_3)}&si Q$\in\{P_1,\dots ,P_{p}
\}$;\cr    k[[t_1,t_2,t_3]]/_{\dis (t_1t_2)}&sinon.\cr}$$}\vs{-2}
\ind Soit $\st_0$ la normalis\'ee de $S_0$ que l'on supposera lisse. Soit
$\varphi\colon\st_0\to\P{}{3}$ le compositum $\st_0\hfl{\nu}{}\st\
\imm{}\P{}{3}$ o\`u
$\nu$ est le morphime de  normalisation. Soit $\lt_0\!=\!\varphi^{-1}L_0$,
le morphisme
$\varphi\!\mid_{\lt_0}\colon\lt_0\to L_0$ est de degr\'e $2$ ramifi\'e le
long de $R_1,\dots, R_p$
qui sont respectivement les images r\'eciproques des points  $P_1,\dots
,P_{p}$.\par
\sec{D\'efinition}{ Soit $X$ une vari\'et\'e lisse de dimension $n$ sur le
corps
$k$. Soient $u_1,\dots ,u_n$ des \'el\'ements de $\Gamma (U,\OO_X)$ o\`u
$U$ est un ouvert de
$X$. On dira que $\{u_1,\dots ,u_n\}$ est un syst\`eme de param\`etres
uniformisants dans $U$ si
$\{du_1,\dots ,du_n\}$ est une base de ${\Omega_{X/k}}\mid_U$.}
\sec{D\'efinition}{Soit $X$ un sch\'ema, et soit $x\in X$ un point. Un
voisinage
\'etale de $x$ dans $X$ est la donn\'ee d'un triplet $(X',e,x')$ o\`u
$e:X'\fh X$ est un morphisme
\'etale de sch\'emas et $x'$ est un point de $X'$ tel que $e(x')= x$.}
\ind Avec ces notations et d\'efinitions on montre la \vs{1}
\sec{Proposition}{ Soit $Q \in L_0 \setminus\{T_1,\dots ,T_t,P_1,\dots
,P_{p}\}$. Soit $F_0$
l'\'equation de
$S_0$; il existe  un voisinage affine \'etale
$(V_0,e_0,Q')$ de $Q$ dans $\P{}{3}$, et il existe $\{ F_0^{\scsc(1)},F_0^{\scsc(2)},
F_0^{\scsc(3)}\}$ un syst\`eme de param\`etres uniformisants autour de
$Q'$ tels que: {\parindent=1cm
\item {\rm (i)} $e_0^*(F_0)\ =\ F_0^{\scsc(1)}. F_0^{\scsc(2)}$.
\item {\rm (ii)} $V(F_0^{\scsc(1)})$ et $V( F_0^{\scsc(2)})$ sont lisses et
se coupent
transversalement en $Q'$.  
\item {\rm(iii)} Quitte \`a se restreindre \`a un ouvert de $V_0$ on a
$L'_0=e_0^{-1}L_0$  est d\'efini par $F_0^{\scsc(1)}=F_0^{\scsc(2)}=0$.}}
\dem Soient $U$ un voisinage affine de $Q$ dans $\P{}{3}$ et
$V=\varphi^{-1}(U)$; notons $R$ un ant\'ec\'edent de
$Q$. Soit $\{u,v\}$ (resp. $\{x,y,z\}$)  un syst\`eme de param\`etres
uniformisants dans $V$ (resp.
dans $U$). On peut supposer (quitte \`a effectuer un changement lin\'eaire
des param\`etres) que
$u(R)=v(R)=0$ et que
$x(P)=y(P)=z(P)=0$. Notons $(\varphi_1,\varphi_2,\varphi_3)$ les
composantes du morphisme $\varphi$.
La matrice du morphisme tangent $T_{\varphi}$ au voisinage de $R$ poss\`ede
un mineur non nul, par
exemple,${\dis\partial{\varphi_1}\over{\dis\partial u
}}{{\dis\varphi_2}\over{\dis\partial v
}}-{\dis\partial{\varphi_1}\over{\dis\partial v
}}{{\dis\varphi_2}\over{\dis\partial u }},$
 alors le morphisme
$$\diagram{V\times\A{}{1}&\hfl{\psi}{}&\hskip -3cm U\cr  (u,v,t)&\mapsto
&(\varphi_1(u,v),\varphi_2(u,v),\varphi_3(u,v)+t)\cr}$$ est \'etale au
voisinage du point $(R,
0)$.\par 
\ind Notons encore $F_0$ la restriction $F_0\mid_U$. Comme l'image de
${\st_0}\times\{0\}$ par $\psi$ est contenue dans $S_0$, on a $\psi^*(F_0)=
t.F_0^{\scsc(1)}$ o\`u
$F_0^{\scsc(1)}$ est une fonction d\'efinie sur $V\times\A{}{1}$. Pour
montrer la proposition, il
suffit de prouver le:\par
\subsec{Lemme}{\it On pose $F_0^{\scsc(1)}= Au + Bv+Ct \mod{(u,v,t)^2}$,
alors
 $ (A,B)\not =(0,0)$.}\par 
\ind En effet si, par exemple, $A\not =0$ et si l'on pose 
$F_0^{\scsc(2)}= t,\ F_0^{\scsc(3)}= v\ {\rm et}\ Q'=(R,0)$,
alors\vs{-2}$$\mid {\partial {(\dis
F_0^{\scsc(1)},F_0^{\scsc(2)},F_0^{\scsc(3)})}\over{\dis\partial (t,u,v)
}}\mid(Q')=\ det\
\pmatrix{C&A&B\cr 1&0&0\cr 0&0&1\cr}=A\not =0$$\vs{-1} et comme
$F_0^{\scsc(1)}(Q')=0$ alors
$\{F_0^{\scsc(1)},F_0^{\scsc(2)},F_0^{\scsc(3)}\}$ est un syst\`eme de
param\`etres uniformisants
autour de $Q'$ (si $B\not =0$ on prend $F_0^{\scsc(3)}= u$). Quitte
\`a restreindre $V$ \`a un voisinage ouvert de $R$, on peut supposer que
$\{F_0^{\scsc(1)},F_0^{\scsc(2)},F_0^{\scsc(3)}\}$ est un syst\`eme de
param\`etres uniformisants
dans $V\times \A{}{1}$.\par
 On pose $V_0=V\times \A{}{1}$ et $e_0$ la compos\'ee de $\psi$ et de
l'immersion
$U\hookrightarrow \P{}{3}$. L'assertion (i) est claire. Pour montrer (ii),
on remarque que
${\dis\partial {\dis F_0^{\scsc(1)}}\over{\dis\partial u }}(Q')=A\not =0$
(i.e.,
$V(F_0^{\scsc(1)})$ est lisse en $Q'$) et que $A$ est un mineur de la
matrice $${\partial {(\dis
F_0^{\scsc(1)},F_0^{\scsc(2)})}\over{\dis\partial (u,v,t)
}}(Q')=\pmatrix{A&B&C\cr 0&0&1\cr}.$$ 
Montrons l'assertion (iii). On a $L_0=V(F_0)$ et $e_0^*(F_0)=
F_0^{\scsc(1)}F_0^{\scsc(2)}$, donc le
lieu des points doubles de $e_0^*(F_0)$ est l'ensemble des z\'eros communs
\`a $F_0^{\scsc(1)}$ et $F_0^{\scsc(2)}$. Donc ensemblistement
$L'_0=V(F_0^{\scsc(1)},F_0^{\scsc(2)})$. Comme $L'_0$ est r\'eduit (car
$L_0$ est un sch\'ema
r\'eduit et $e_{\scsc0}$ est \'etale) et
$V(F_0^{\scsc(1)},F_0^{\scsc(2)})$ est r\'eduit (car $F_0^{\scsc(1)}$ et
$F_0^{\scsc(2)}$ font partie d'un syst\`eme de param\`etres uniformisants),
alors l'\'egalit\'e
$L'_0=V(F_0^{\scsc(1)},F_0^{\scsc(2)})$ est une \'egalit\'e sch\'ematique.
\smallskip 
{\it Preuve du lemme 2.4.1. }On va montrer que
$A\hs{-1}=\hs{-1}B\hs{-1}=\hs{-1}0$ implique
$\mid {\dis\partial {(\dis e^*_{\scsc 0}x,e^*_{\scsc 0}y,e^*_{\scsc
0}z)}\over{\dis\partial (u,v,t)
}}\mid (Q')=\hs{-1}0$, ce qui contredit le fait que
$e_0$ est \'etale en $Q'$. Posons
$$\system {x^*(u,v,t)=e^*_{\scsc 0}x=a_ {\scsc 100}u\ +\ a_ {\scsc 010}v\
+\ a_ {\scsc 001}t\ +\
{\rm\ des\ termes\ de}\  (u,v,t)^2\cr  y^*(u,v,t)=e^*_{\scsc
0}y=b_{\scsc100}u\ +\ b_{\scsc010}v\
+\  b_{\scsc001}t\ +\ {\rm des\ termes\ de}\  (u,v,t)^2\cr}$$  On a alors $$\mid {\dis\partial
{(\dis e^*_{\scsc 0}x,e^*_{\scsc 0}y,e^*_{\scsc 0}z)}\over{\dis\partial
(u,v,t) }}\mid
(Q')=\alpha\cdot{\dis\partial{\dis z^*}\over{\dis\partial u }}
(Q')+\beta\cdot{\dis\partial{\dis
z^*}\over{\dis\partial v}} (Q')+\gamma\cdot{\dis\partial{\dis
z^*}\over{\dis\partial t}}(Q')$$ o\`u
$\alpha=a_{\scsc 010}b_{\scsc001}-a_{\scsc 001}b_{\scsc010}$,
$\beta=a_{\scsc 100}b_{\scsc
001}-a_{\scsc 001}b_{\scsc100}$ et $\gamma=a_{\scsc 100}b_{\scsc
010}-a_{\scsc 010}b_{\scsc100}$.\par
\ind D'autre part on a:\par
$\ x^*.y^*=a_{\scsc 100} b_{\scsc100}u^2+ a_{\scsc 010}b_{\scsc010}v^2+ (a_
{\scsc 010}
b_{\scsc100}+a_{\scsc 100} b_{\scsc010})uv+{\ \rm des\ termes\ de}\ 
(u,v)^3 +\
 t[(a_{\scsc 001} b_{\scsc100}+a_ {\scsc 100} b_{\scsc001})u\ +\ (a_{\scsc
001}
b_{\scsc010}+a_{\scsc 010} b_{\scsc001})v\ +\ a_{\scsc 001} b_{\scsc 001}t\
+\ {\rm des\ termes\
de}\  (u,v,t)^2]$. \par
\abovedisplayskip=4pt plus 3pt minus 9pt
\abovedisplayshortskip=-8pt plus 3pt 
\belowdisplayskip=12pt plus 3pt minus 9pt
\belowdisplayshortskip=9pt plus 3pt minus 4pt
 \ind Quitte \`a changer les coordonn\'ees, on sait que $\hat\OO_{S_0,
Q}\simeq
k[[x,y,z]]/_{\displaystyle (xy)}$ autrement dit, le d\'eveloppement de
Taylor de $F_0$ est \'egal \`a
$xy$, ainsi $x^*y^*=$$e_0^*F_0= t.F_0^{\scsc(1)}$, ce qui implique en
particulier  les \'egalit\'es
suivantes:\par
\hs{6} $ A=a_ {\scsc 001} b_{\scsc100}+a_ {\scsc 100} b_{\scsc001},\ 
\ B=a_{\scsc 001} b_{\scsc010}+a_{\scsc 010} b_{\scsc001},\ \  a_{\scsc
100}b_{\scsc100}= a_ {\scsc
010}b_{\scsc010}= a_ {\scsc 010}b_{\scsc100}+a_ {\scsc
100}b_{\scsc010}=0.$\par
 Celles-ci donnent dans le cas $A=B=0$ les \'egalit\'es \smallskip
\ind $ \alpha=-2a_ {\scsc 001} b_{\scsc010}=2a_ {\scsc 010} b_{\scsc001},\ 
\  \beta=-2a_ {\scsc
001} b_{\scsc100}=2a_ {\scsc 100}b_{\scsc001},\ \  
\gamma=-2a_ {\scsc 010}b_{\scsc100}=2a_ {\scsc 100}b_{\scsc010}$\par d'o\`u
l'on d\'eduit
$\alpha=\beta=\gamma=0$ en \'etudiant les 4 cas donn\'es par les
\'egalit\'es
$a_{\scsc 100}b_{\scsc100}= a_ {\scsc 010}b_{\scsc010}=0$.\cqfd\vs{5}
\ind La premi\`ere question que l'on pourrait se poser \`a propos de la
proposition pr\'ec\'edente
est de savoir si celle-ci est est valable pour les d\'eformations plates 
infinit\'esimales de
$S_0$; le cas \'ech\'eant, la surface $S_0$ ne pourra pas \^etre
lissifi\'e. Nous allons montrer que
la r\'eponse est ``oui" sous certaines conditions relatives  \`a
l'existence de courbes trac\'ees
sur $S_0$ d\'ecoupant sur $L_0$ un nombre ``assez suffisant" de points qui
sera au moins le degr\'e
de la surface duale \`a $S_0$. L'utilisation des sch\'emas polaires
d\'efinis dans la section 1 est
la base sur laquelle reposera toute la d\'emonstration. Nous aurons au
pr\'ealable besoin de\par   
\sec{Notations et hypoth\`eses suppl\'ementaires}{} 
 Soit $S_0$ comme dans (2.1), on supposera de plus que
$L_0$ est lisse et irr\'eductible. (Supposer que $L_0$ est irr\'eductible
ne nuira pas \`a la
g\'en\'eralit\'e dans la mesure o\`u on pourra remplacer dans ce qui suit
$L_0$ par une de ses
composantes irr\'eductibles.)\par
  \ind Soit $C_0\subset S_0$ une courbe lisse et irr\'eductible qui coupe
$L_0$ en dehors des points-pinces et ce transversalement, on note $\{
Q_1,\dots , Q_r\}= C_0\cap
L_0$. Les points $Q_1,\dots,Q_r$ n'\'etant pas des points-pinces, on peut
leur appliquer la
proposition (2.4). Il existe donc un voisinage \'etale affine
$(V_0,e_0,Q'_1,\dots,Q'_r)$ des points $Q_1,\dots,Q_r$ dans $\P{}{3}$, avec
$Q_i=e_0(Q'_i)$ pour tout $i$, et des  param\`etres uniformisants $\{
F_0^{\scsc(1)},
F_0^{\scsc(2)},\ F_0^{\scsc(3)}\}$ dans $V_0$ tels que: {\parindent=1cm
\item {\rm (i)} $e_0^*(F_0)\ =\ F_0^{\scsc(1)}. F_0^{\scsc(2)}$.
\item {\rm (ii)} $V(F_0^{\scsc(1)})$ et $V( F_0^{\scsc(2)})$ sont lisses et
se coupent
transversalement en $Q'$.  
\item {\rm(iii)} $L'_0=e_0^{-1}L_0$ est d\'efini par
$F_0^{\scsc(1)}=F_0^{\scsc(2)}=0$.}\par Pour tout $n\geq0$, on pose
$V_n=\fibp{V_0}{}{D_n}$ et
$e_n=e_0\times1_{D_n}$, donc
$e_n\colon V_n\to\pnd{n}$ est un rev\^etement \'etale de $\pnd{n}$, et
puique $V_n$ (resp.$\pnd{n}$)
a le m\^eme espace sous-jacent que $V_0$ (resp. $\P{}{3}$) on peut
consid\'erer $Q'_1,\dots,Q'_r$
(resp. $Q_1,\dots,Q_r$) comme des points de $V_n$ (resp.
$\pnd{n}$), ainsi $(V_n,e_n,Q'_1,\dots,Q'_r)$ est un voisinage \'etale dans
$\pnd{3}$ des points
$Q_1,\dots,Q_r$. De plus, pour tout $n\geq0$, on a
$V_{n+1}=\fibp{V_n}{D_{n+1}}{D_n}$ et
$e_{n+1}\mid_{V_n}=e_n$.\par
\ind On note $\eps{n}{}$ l'image de $t$ dans $k[t]/_{\dis(t^{n+1})}$\par
\sec{Proposition}{ On garde les notations et les hypoth\`eses
pr\'ec\'edentes. On suppose qu'il
existe une famille de drapeaux
$(C_n\subset S_n\subset\pnd{n})_{n\geq0}$, d\'eformations infinit\'esimales
plates du drapeau
$(C_0\subset S_0\subset\P{}{3})$ avec, pour tout $ n\geq 0$,
$S_{n}=\fibp{S_{n+1}}{D_{n+1}}{D_{n}}$
et $C_{n}=\fibp{C_{n+1}}{D_{n+1}}{D_{n}}
$. Notons $s$ le degr\'e de $S_0$ et, pour tout $n\geq0$, 
$F_n\in\coh{0}{\pnd{n}}{\OO_{\pnd{n}}(s)}$ l'\'equation de $S_n$ dans
$\pnd{n}$.\par 
 \ind Si $r=\#(C_0\cap L_0)>s(s-1)^2$, alors pour tout
$n\geq 0$, il existe des fonctions $F_{n}^{\scsc(1)}$ et $F_n^{\scsc(2)}$
d\'efinies sur
$V_n$ telles que:\par
\ind {\rm (i)} $F_0^{\scsc(1)}$ et $F_0^{\scsc(2)}$ sont les param\`etres
d\'efinis dans (2.4)\par
\ind {\rm (ii)} $e^*_{\scsc n}F_n\ =\ F_{n}^{\scsc(1)} F_n^{\scsc(2)}$.\par
\ind {\rm (iii)} $F_{n}^{\scsc(1)}\mid_{V_{n-1}}=F_{n-1}^{\scsc(1)}$ et
$F_n^{\scsc(2)}\mid_{V_{n-1}}= F_{n-1}^{\scsc(2)}$.\par
\ind {\rm (iv)} $F_{n}^{\scsc(1)}$ et $F_n^{\scsc(2)}$ font partie d'un
syst\`eme de param\`etres
uniformisants dans $V_n$.} 
\dem On fait une r\'eccurence sur $n$. Pour $n=0$, on applique la
proposition (2.4).\par
 Soit $n\geq 1$, on suppose le r\'esultat vrai jusqu'\`a l'ordre $n\!-\!1$.
Alors en particulier, il
existe deux fonctions $F_{n-1}^{\scsc(1)}$ et
$F_{n-1}^{\scsc(2)}$ sur $ V_{n\!-\!1}$ faisant partie d'un syst\`eme de
param\`etres uniformisants
telles que\hfill\break\vs{-3}
\hs{30}$e^*_{\scsc n-1}F_{n-1} =F_{n-1}^{\scsc(1)} F_{n-1}^{\scsc(2)}$ 
$F_{n}^{\scsc(1)}\mid_{V_{n-1}}=F_{n-1}^{\scsc(1)}$ et
$F_n^{\scsc(2)}\mid_{V_{n-1}}= F_{n-1}^{\scsc(2)}$.\smallskip 
 Par construction, $V_{n-1}$ est d\'efini comme sous-sch\'ema ferm\'e du
sch\'ema affine
$V_n$ par
$\eps{n}{n}=0$. On peut alors \'etendre arbitrairement les fonctions
$F_{n-1}^{\scsc(1)}$ et
$F_{n-1}^{\scsc(2)}$ en deux fonctions $u$ et $v$ d\'efinies sur $V_n$. On
a $( e_{\scsc n}^*F_n -
uv)\mid_{V_{n-1}}=0$, donc il existe une fonction $h$ d\'efinie sur
$V_n$ telle que 
$$e_{\scsc n}^*F_n = uv+\eps{n}{n} h\leqno (2.5.1)$$  Et la proposition est
d\'emontr\'ee si on
prouve le:\smallskip
\sec{Lemme}{ La fonction $h$ d\'efinie dans (2.5.1) est \'egale \`a $0$
modulo l'id\'eal
$(u,v,\eps{n}{})$.}
\hs{4} En effet, si $h=au+bv\ [\eps{n}{}]$ alors $ e_{\scsc
n}^*F_n\hs{-1}=\hs{-1}uv+\eps{n}{n}au
+\eps{n}{n}b v$. On pose
$F_{n}^{\scsc(1)}=\hs{-1}u+\eps{n}{n}b$ et
$F_{n}^{\scsc (2)}=v+\eps{n}{n}a$,\hfill\break alors $e_{\scsc
n}^*F_n=F_{n}^{\scsc(1)}F_{n}^{\scsc
(2)}$, d'o\`u les assertions (i) et (ii). L'assertion (iii) est \'evidente.
Pour voir l'assertion
(iv),  on a
$\tens{\Omega_{\scsc V_n/D_n}}{B_n}{k}\simeq\Omega_{\scsc V_0/L_0}$, comme
$du\mid_{L_0}(=dF_0^{\scsc(1)})$ et $dv\mid_{L_0}(=dF_0^{\scsc(2)})$ sont
lin\'eairement
ind\'ependants dans $\Omega_{\scsc V_0/L_0}$ alors, par le lemme de
Nakayama, $du$ et
$dv$ le sont dans $\Omega_{\scsc V_n/D_n}$, donc aussi $dF_n^{\scsc(1)}$ et
$dF_n^{\scsc(2)}$ car
$${\partial
\dis F_{n}^{\scsc(1)}\over{\partial \dis u }}{\partial \dis
F_n^{\scsc(2)}\over{\partial \dis v
}}-{\partial \dis F_{n}^{\scsc (1)}\over{\partial
\dis v }}{\partial
\dis F_n^{\scsc(2)}\over{\partial \dis u }}=1+\eps{n}{n}({\partial \dis
b\over{\partial\dis u
}}+{\partial \dis a\over{\partial \dis v }})$$  est une unit\'e.\par  {(\it
Preuve du lemme)} On
proc\`edera en trois \'etapes: \par     {$\underline{\bf 1^{\b o}\acute
Etape}$:} On va montrer que
pour tout
$i\in\{1,\dots,r\}$ on a $h(Q'_i)=0$.\smallskip  
\ind Pour tout $j\geq 0$ notons $C'_j$ l'image inverse par $e_{\scsc j}$ de $C_j$, i.e.,
$C'_j=\fibp{C_j}{\pnd{j}}{V_j}$. Fixons un point $Q$ quelconque parmi les
$Q_i$, il suffit de
montrer $h\mid_{C'_0}(Q')=0$.\par On a
$F_{n-1}\mid_{C_0}=0$ et $e_{\scsc
n-1}^*F_{n-1}=F_{n-1}^{\scsc(1)}F_{n-1}^{\scsc(2)}$. Comme
$\OO_{{C'_0},{Q'}}$ est int\`egre ($C'_0$ est lisse, car \'etale sur $C_0$
qui est lisse par
hypoth\`ese), alors dans $\OO_{{C'_0},{Q'}}$ on a $
F_{n-1}^{\scsc(1)}\mid_{C'_0}=F_0^{\scsc(1)}\mid_{C'_0}=0$ ou
$F_{n-1}^{\scsc(2)}\mid_{C'_0}=F_0^{\scsc(2)}\mid_{C'_0}=0$. On suppose,
par exemple,
$F_{n-1}^{\scsc(1)}\mid_{C'_0}=0$. Les courbes $C_0$ et $L_0$ sont
transverses en
$Q$, donc $C'_0$ et $L'_0=e_0^{-1}(L_0)$ le sont aussi en $Q'$, et comme
$L'_0$ est d\'efinie par
l'id\'eal $(F_0^{\scsc(1)},F_0^{\scsc(2)})$, alors
$F_{n-1}^{\scsc(2)}\mid_{C'_0}\not =0$ dans $\OO_{{C'_0},{Q'}}$. Par
ailleurs, dans
$\OO_{{C'_{n-1}},{Q'}}$ on a $F_{n-1}^{\scsc(1)}
F_{n-1}^{\scsc(2)}\mid_{C'_{n-1}}=0$ (car
$C_{n-1}\subset S_{n-1}$), et $F_{n-1}^{\scsc(2)}\mid_{C'_0}\not =0$
implique
$F_{n-1}^{\scsc(1)}\mid_{C'_{n-1}}=0$; pour le voir on utilisera le fait
que
$\OO_{C'_{n-1},Q'}$ est isomorphe \`a $\tens{\OO_{C'_0,Q'}}{k}{k[t]/_{\dis
(t^n)}}$ puisque
localement les d\'eformations infinit\'esimales des courbes lisses sont
triviales. Dans
$\OO_{C'_n,Q'}$, via un isomorphisme
$\OO_{C'_n,Q'}\simeq\tens{\OO_{C'_0,Q'}}{k}{k[t]/_{\dis (t^{n+1})}}$
choisi, on
\'ecrit:\par 
$u\mid_{C'_n}=\sum_{i=1}^{n}\eps{n}{i}u_i,\  v\mid_{C'_n}=
\sum_{i=1}^{n}\eps{n}{i}v_i \ {\rm  et}\
h\mid_{C'_n}=\sum_{i=1}^{n}\eps{n}{i}h_i$ o\`u
$u_i,\ v_i$ et $h_i$ sont dans $\OO_{C'_0,Q'}$ pour tout $i=1,\dots,n$. On
a alors
$F_{n-1}^{\scsc(1)}\hs{-1}\mid_{C'_{n-1}}\hs{-1}=\hs{-1}u\mid_{C'_{n-1}}=
\sum_{i=1}^{n-1}\eps{n-1}{i}u_i$, et
$F_{n-1}^{\scsc(1)}\mid_{C'_{n-1}}=\hs{-1}0$ implique
$u\mid_{C'_n}= u_n \eps{n}{n}$, et comme $e^*_{\scsc n} F_n
\mid_{C'_n}=( uv+\varepsilon_{n}^n h)\mid_{C'_n}=0$
($C_n\hs{-1}\subset\hs{-1} S_n$), alors
$\varepsilon_{n}^n(u_nv_0 + h_0) \hs{-1}=\hs{-1}0$ dans $\OO_{C'_n,Q'}$,
autrement dit,
$u_nv_0\hs{-.5}+\hs{-.5}h_0\hs{-1}=\hs{-1}0$ dans $\OO_{C'_0,Q'}.$
Rappelons que la courbe $L'_0$
est d\'efinie par l'i\'eal $(F_{0}^{\scsc(1)}, F_{0}^{\scsc(2)})$, donc
$v_0(Q')=F_{0}^{\scsc(2)}\mid_{C'_0}(Q')= 0 $, d'o\`u
$h_0(Q')=-(u_nv_0)(Q')=0.$\par

{$\underline{\bf 2^{\b o}
\acute Etape}$:} On va se placer dans la situation de la  proposition
(1.10). Soit $\Delta$ une
droite g\'en\'erale de $\P{}{3}$, on consid\`ere
$\eta_0\colon\wp\to\P{}{3}$ l'\'eclatement de
$\P{}{3}$ en
$\Delta$ et $q_0\colon\wp\to\P{}{1}$ la projection induite sur $\P{}{1}$.
\par
\ind Soit $G\in \coh{0}{\P{}{3}}{\OO_{\P{}{3}}(s)}$ tel que la surface
$V(G)$ soit lisse. Soit
$F\in\cohr{0}{\OO_{\P{}{3}\times\A{}{1}}(s)}$ tel $G=F\mid_{\P{}{3}\times
\{1\}}$ et
$F_n=F\mid_{\P{}{3}\times D_n}$.( Le polyn\^ome $F$ existe car
l'application
naturelle$$\cohr{0}{\OO_{{\P{}{3}}\times\A{}{1}}(s)}\hfl{\Phi}{}
\cohr{0}{\OO_{\P{D_n}{3}}(s)}\times\cohr{0}{\OO_{\P{}{3}\times\{1\}}(s)}$$
est surjective par le
lemme chinois.) \par
 On pose $\wpol=\pol(q_0\times id_{\A{}{1}},(\eta_0\times
id_{\A{}{1}})^*F)$, et
$\pol_n=\pol(q_n,\eta_n^*F_n)$ avec $\eta_n=\eta_0\times id_{D_n}$ et
$q_n=q_0\times id_{D_n}$. Nous
avons alors le\par  
\subsec{Sous-lemme}{ On suppose $h\not\equiv 0 \pmod{(u,v,\eps{n}{})}$.
Soit $R'\in L'_0$ tel que
$h(R')=0$. Notons $e_n(R')=R$, alors \par 1) Il existe un ouvert $\UU
\subset \pol_n$ tel que  $\dim
(\pol_n\setminus \UU )=0$ autour de $R$ et $\UU\subset \wpnd{n-1}$.\par 2)
$R$ est un point
immerg\'e de $\pol_n$; plus pr\'ecis\'ement, le support de la fonction
$\varphi:=\eps{n}{n}\mid_{\pol_n}$ contient $R$ comme point isol\'e.\par 3)
Soit $\w\varphi=
t^n\mid_{\wpol}$, alors Supp ($\w\varphi$)$\cap (\wp\times \{0\})$ contient
$R\times\{0\}$ comme
point isol\'e.} 
 \dem La droite $\Delta$ peut
\^etre choisie de telle mani\`ere qu'elle soit disjointe de $L_0$. Ainsi le
point $R$ peut \^etre vu
comme point de $\pol$ (car $\P{}{3}\setminus \Delta\simeq \wp \setminus
\w\Delta$), et $R$ est un point de $\pol_n$ car $R\in L_0$ (cf. 1.8).\par 
(1) On se place dans un
voisinage de $R$. Comme on est en dehors de $\Delta$, alors le morphisme
$q_0$ au voisinage de $R$
n'est autre que la projection de centre
$\Delta$, i.e, $q_0(x,y,z)=z$ o\`u
$(x,y,z)$ est un syst\`eme de param\`etres locaux autour de $R$.
Localement,
$\pol_0$ est d\'efini par l'id\'eal $(F_0,{ \dis\partial
F_0\over\dis\partial x},{
\dis\partial F_0\over\dis\partial y})$, donc dans $V_n$ le sch\'ema
$\pol_n'$, pr\'eimage de
$\pol_n$, est  d\'efini autour de $R'$ par l'id\'eal \smallskip
\hs{50}$(uv+ \eps{n}{n} h,v+\eps{n}{n}{
\dis\partial h\over\dis\partial u}, u+\eps{n}{n}{ \dis\partial
h\over\dis\partial v})$.\smallskip Ce
dernier \'etant contenu dans $(u,v,\eps{n}{})$, on a
$h\not=0\pmod{(u,v,\eps{n}{})}$ sur $\pol_n'$. Posons\vs{2}  
\hs{2}$\UU = \{ M\in \pol_n$ tel qu'il existe $ M'\in \pol'_n$  au-dessus
de
$M$ avec $ h(M')\not =0\pmod{(u,v,\eps{n}{})}\}$.\vs{2}  Consid\'erons
$\UU'=e_n^{-1}\UU$;
ensemblistement $\UU'=\{M'\in \pol_n'\mid\ h(M')\not
=0\pmod{(u,v,\eps{n}{})}\}$. Alors $\UU'$ est
un ouvert non vide de
$\pol_n'$ (par hypoth\`ese), ainsi $\UU$ est un ouvert non vide de $\pol_n$
et donc dim
$(\pol_n\setminus \UU )= 0$. Soit $M'\!\in \UU'$, on a par
d\'efinition\smallskip
\hs{38}$-\eps{n}{n}h(M')=uv(M')= (-\eps{n}{n}{\dis{\partial h\over\partial
v}})(-\eps{n}{n}{\dis{\partial h\over\partial u}})(M')=0$\smallskip
 et puisque la fonction $h$ est inversible sur $\UU'$, on a
$\eps{n}{n}\mid_{\UU'}=0$, et donc
$\eps{n}{n}\mid_{\UU}=0$. Mais
$\eps{n}{n}=0$ est l'\'equation de $D_{n-1}$ dans $D_n$, alors
$\UU\subset\wpnd{n-1}$.\par (2) Il
suffit de montrer que (a) dim (Supp $\varphi)=0$ et que (b) $\varphi\not
=0$ au voisinage de
$R$.\par
  \ind Pour voir (a) consid\'erons un point $x\in \UU$, alors
$\varphi=0$ au voisinage de $x$ puisque, d'apr\`es (1) ci-dessus, $\UU
\subset \pnd{n-1}$, donc
$x\not\in$ Supp $\varphi$. D'o\`u Supp $\varphi\subset Pol\setminus \UU$ et
est donc de dimension
$0$.\par
\ind Montrons (b). Comme tous les calculs sont locaux, on confondra $R$ et
$\pol_n$ respectivement
avec $R'$ et $\pol_n'$. Soit
$(u',v',z)$ le syst\`eme de coordonn\'ees locales autour de $R$ donn\'e par
\smallskip \hs{50}$u'=u+\varepsilon_{n}^n{\dis{\partial h\over\partial
v}},\  v'=v+\varepsilon_{n}^n
{\dis{\partial h\over\partial u}}$\smallskip 
 Soit $\Sigma$ le sch\'ema d\'efini par $u'\!\!=\!\! v'
\!\!=\!\!z\!\!=\!\!0$. Alors $\Sigma\simeq
D_n$, (ensemblistement, c'est $R$ et sch\'ematiquement c'est le point
immerg\'e), et
$\varphi\mid_{\Sigma} \not =0$. Il reste \`a voir que $\Sigma$ est un
sous-sch\'ema ferm\'e de
$\pol_n$. Puisque $h(R)\!=\! 0$, on a
$h\! =\! au'+bv'+ cz+d\eps{n}{}\ $, et donc sur $\Sigma$ on a
$uv+\eps{n}{n}{\dis{\partial
h\over\partial v}}=\varepsilon_{n}^n(au'+bv'+cz+d\varepsilon_{n})$;
autrement dit, $\pol_n \cap
\Sigma=\Sigma$. D'o\`u l'immersion $\Sigma\subset\pol_n$.\vs{1}  

 (3) On se place dans un voisinage de $R$. On a $\pol_n =\wpol\cap
({\w{\P{}{3}}}\times D_n) $
(cf. 1.10), donc $\OO_{\pol_n, {R}}=\OO_{\wpol, {R}}/_{\dis(t^{n+1})}$ et 
$\varphi =
\w\varphi\pmod{t^{n+1}}$ ( car $\w\varphi\mid_{\pol_n}=\varphi$, ainsi
$\w\varphi \not =0$ (car
$\varphi \not =0$) et $R\in $ Supp $\w\varphi$ (car $R \in $ Supp
$\varphi$). Reste \`a voir que $R$ est un point isol\'e dans Supp
($\w\varphi$)$\cap (\wp\times
\{0\})$. Pour cela consid\'erons $x\in \wpol$ un point g\'en\'eral, par
exemple
$x\in\UU\cap\wpol$ o\`u $\UU$ est l'ouvert d\'efini dans le paragraphe (1)
ci-dessus. Montrons que
$x\not\in$ Supp
$\w\varphi$.\par\ind On a $x\not =R$ (car $R\not\in\UU$), et puisque
$\varphi\mid_{\UU}=0$ alors
$\varphi_x =0$. Il suffit donc de montrer que le morphisme 
$\OO_{\wpol, x}\to\OO_{\pol_n,x}=\OO_{\wpol,x}/_{\dis(t^{n+1})}$ est
injectif (car il envoie
$\w{\varphi_x}$ sur $\varphi_x$). \par
\ind On a $\OO_{\pol_n,x} = \OO_{\wpol , x}/_{\dis{t^{n+1}.\OO_{\wpol,
x}}}$, or $\UU\subset
{\w\P{}{3}}\times D_{n-1}$, donc $t^n .\OO_{\pol_n,x}=0 $, et donc
$t^n.\OO_{\wpol , x}\subset t^{n+1}.\OO_{\wpol , x}$. Ceci entra\^ine que pour tout $ m\geq n$ on a
$t^m.\OO_{\wpol , x}\subset t^{m+1}.\OO_{\wpol , x}$, i.e., \par
 $ t^n.\OO_{\wpol , x}\subset\bigcap_{m\geq 1} t^{m+n}.\OO_{\wpol , x}$; et
par le th\'eor\`eme de
Krull, ceci implique $ t^n.\OO_{\wpol , x}=0$. D'o\`u $\OO_{\pol_n
,x}=\OO_{\wpol,x}$.
\vs{2} {$\underline{\bf 3^{\b o}\acute Etape\ et\ fin}$:} On suppose
$h\not\equiv 0
\pmod{(u,v,\eps{n}{})}$. Alors d'apr\`es (2.7.1, (3)) et (1.10, (v)), le
nombre des z\'eros de
$h\mid_{L'_0}$ est au plus $s(s-1)^2$. Ainsi l'hypoth\`ese $r>s(s-1)^2$ et
le r\'esultat de la
1\`ere \'etape impliquent le lemme.\cqfd\vs{4}
\ind Ce que dit, plus ou moins implicitement, la proposition pr\'ec\'edente
est que la condition
$r>s(s-1)^2$ suffit pour qu'aucune d\'eformation plate de $S_0$ contenant
une d\'eformation plate 
de $C_0$ ne puisse lissifier le lieu double de $S_0$. Plus pr\'ecis\'ement
on montre la \vs{2}
\sec{Proposition}{Soit $S_0\subset\P{}{3}$ une surface \`a singularit\'es
ordinaires comme dans
(2.5), dont la courbe double $L_0$ est lisse et connexe. Soit $C_0\subset
S_0$ une courbe lisse et
connexe qui coupe $L_0$ en dehors des points-pinces en $r$ points et ce
transversalement. On suppose
qu'il existe une d\'eformation plate $\ee C$ (resp. $\ee S$) de $C_0$
(resp. $S_0$) dans
$\P{}{3}$ via $A=k[[t]]$, avec ${\ee C}\subset \ee S$. Notons $s$ le
degr\'e de $S_0$. Si
$r>s(s-1)^2$, alors
$\ee S$ est singuli\`ere le long d'un sous-sch\'ema (ferm\'e) $\ee L$
avec\par
\hs{7}{\rm (i)} $\ee L$ est une d\'eformation plate de $L_0$ sur $\ZZ=\spec
A$.\par
\hs{7}{\rm (ii)} La fibre g\'en\'erique de
$\ee S$ est une surface singuli\`ere contenant la fibre g\'en\'erique de
$\ee L$ comme courbe
double.}
\dem Soit ${\ee L}\subset {\ee S}$ le lieu o\`u la projection ${\ee S}\fh
\ZZ$ n'est pas lisse,
c'est un sous-sch\'ema ferm\'e de ${\ee S}$ tel que ${\ee L}\cap
(\P{}{3}\times D_0)=L_0$ et ${\ee
L}\cap\ (\P{}{3}\times D_m)=L_m$ est le lieu singulier de la projection
$S_m\fh D_m$ o\`u $S_m= {\ee
S}\mid_{\P{D_m}{3}}$.\par (i)\hs{8} $Pas\ 1.$ Pour tout entier naturel $m$,
le morphisme $L_m\to
D_m$ est lisse en dehors des points-pinces. En effet, nous sommes dans les
hypoth\`eses de la
proposition (2.6), en particulier $r>s(s-1)^2$, donc, en dehors des
points-pinces, $L_m$
est d\'efinie par
$u_m$$=$$v_m$$=$$0$$\mod{\eps{m}{}}$ \`a un rev\^etement \'etale pr\`es. Or
$u_0$ et $v_0$ sont
transverses, donc Jac $(u_0,v_0)\not =0$, ainsi Jac
$(u_m,v_m)\not =0\mod{\eps{m}{}}$, et ceci implique que l'intersection des
nappes $(u_m$$=$$0)$ et
$(v_m$$=$$0)$ est lisse sur $D_m$.\par
\ind $Pas\ 2.$ La courbe $L_0$ n'est pas ensemblistement une composante
irr\'eductible de $\ee L$.
Sinon, consid\'erons $U$, une composante irr\'eductible de $\ee L$ telle
que
$U_{\rm red}=L_0$. Soit $\NN$ l'id\'eal de $L_0$ dans $U$, (i.e. le
nilradical de
$\Gamma(U,\OO_{U}$). Alors il existe un entier$n_0$ tel que
${\NN}^{n_0}=0$.  Par ailleurs, pour
tout $m\geq 0$, on a $L_m\simeq L_0\times D_m$ (dans la topologie \'etale)
car
$L_m\fh D_m$ est lisse en dehors des points-pinces (cf. le pas 1), donc $
(L_m)_{red}= L_0$. Ainsi 
$L_m\subset U $. On va voir que pour $m>n_{0}$, ceci est impossible:  on a
$t\mid_{L_0}=0$ (car
$\P{k}{3}$ est le sous-sch\'ema ferm\'e de $\P{A}{3}$ d\'efini par
l'\'equation $t=0$), ce qui
implique $t\mid_U\in{\NN}$, donc $t^{n_0}\mid_U=0$, et, comme on vient de
le montrer,
$t^{n_0}\mid_{L_m}=0$ pour tout $m$, mais ceci est impossible si $m\geq
{n_0}$ car
$t^{n_0}\mid_{D_m}= \eps{m}{n_0}\not =0$.\par
\ind $Conclusion.$ Le morphisme ${\ee L}_{red}\fh \ZZ$ est plat car la
fibre au-dessus du point
ferm\'e, \`a savoir la droite $L_0$, n'est pas ensemblistement une
composante irr\'eductible de
${\ee L}$ (cf. le pas 2). Donc le sch\'ema  $\LL$ est plat sur  $\ZZ$. Ceci
d\'emontre l'assertion
(i).\vs{2} (ii)\hs{3} Notons $b$ le point ouvert de $\ZZ$ et ${\ee L}_b$ la
fibre (g\'en\'erique) de
$\ee L$ au-dessus de $b$. Montrons que le lieu double de ${\ee S}_b$ est
ensemblistement \'egal \`a
${\ee L}_b$.\par
\ind On a ${\ee L}_b\buildrel {\scsc ensembl.} \over =L_0$. Notons $\JJ =
\JJ_{L_0/\P{}{3}}$ et $f_b$ l'\'equation de ${\ee S}_b$. Puisque $L_0$ est
singuli\`ere dans
${\ee S}_b$, alors
$f_b \in H^0(\JJ^2 (s))$. Soit $\bar f_b$$\in \cohr{0}{\JJ^2/_{\dis\JJ^3}
(s)}$ =
$\cohr{0}{{\rm Sym}^2\ \NN_{L_0/\P{}{3}}^{\vee}(s)}$ l'image de $f_b$.
Alors $\bar f_b$ est une
forme quadratique. Dire que $L_0$ est une courbe double, \'equivaut \`a
dire que la forme
quadratique $\bar f_b$ est  non d\'eg\'en\'er\'ee au  point g\'en\'eral de
$L_0$, i.e., de
discriminant non nul. Si c'est vrai au point point ferm\'e de $\ZZ$, c'est
vrai au point
g\'en\'erique $b$, par semi-continuit\'e du discriminant.\cqfd\vs{4}
\ind Dans la proposition pr\'ec\'edente, l'existence de la surface $\ee S$,
qui est une hypoth\`ese
essentielle, est en g\'en\'eral tr\`es difficile \`a r\'ealiser. Nous
allons voir dans ce qui suit,
qu'en imopsant certaines conditions sur le genre et le degr\'e de
$C_0$, cette hypoth\`ese est toujours vraie. Plus pr\'ecis\'ement nous
avons la\vs{2}
 \sec{Proposition}{Soient $S_0$, $L_0$, $C_0$, $r$, $s$ comme dans la
proposition pr\'ec\'edente. On
note $d$ et $g$ respectivement le degr\'e et le genre de $C_0$. On suppose 
$r>s(s-1)^2$. Si $d>s^2$
et $g>G(d,s+1)$ alors, pour toute d\'eformation plate $\ee C$ dans
$\P{}{3}$ de $C_0$ via $A=k[[t]]$, il existe une surface ${\ee
S}\subset\P{A}{3}$ de degr\'e $s$
telle que\par  (i) $\ee S$ est une d\'eformation plate de $S_0$ via $A$
contenant $\ee C$ comme
sous-sch\'ema ferm\'e.\par   (ii) $\ee S$ est une surface singuli\`ere de
lieu singulier une courbe
double, d\'eformation plate de $L_0$.
\par}
\dem La conclusion de (i) est une cons\'equence imm\'ediate du lemme (2.10)
ci-dessous  appliqu\'e
\`a
$n=s$. On en d\'eduit (ii) d'apr\`es (2.8).\cqfd\vs{2}\sec{Lemme}{Soit
$C\subset \P{}{3}$ une courbe
lisse et connexe de degr\'e
$d$ et genre $g$, trac\'ee sur une surface int\`egre $S$ de degr\'e $s$. On
suppose
$g>G(d,n+1)$ pour un certain entier $n\geq s$. Alors, pour toute
d\'eformation plate $\ee C$ de
$C$ dans $\P{}{3}$, param\'etr\'ee par un anneau de valuation discr\`ete
$A$ sur
$k$, il existe une famille $\ee T$ de surfaces de $\P{A}{3}$, de degr\'e
$n$ telle que:
\decale{\rm (i)} ${\ee C}$ est un sous-sch\'ema ferm\'e de $\ee T$ plat sur
$\ZZ=\spec A$.
\decale{\rm (ii)} Si de plus $d\hs{-1}>\hs{-1}s.n$, alors la fibre
sp\'eciale de $ \ee T$ est
\'egale \`a $S\cup\Lambda$, o\`u $\Lambda$ est une surface de degr\'e
$n-s$.}
\dem  Pr\'ecisons d'abord quelques notations: on note $a$, (resp. b) le
point ferm\'e (resp. le
point ouvert) de $\ZZ$. Soit $\pi$ un g\'en\'erateur de l'id\'eal maximal
de $A$ tel que $v(\pi)=1$
o\`u $v$ d\'esigne la valuation de $A$. Soit $K$ le corps de fractions de
$A$.\par
\ind La fibre g\'en\'erique ${\ee C}_b=\fibp{{\ee C}}{\ZZ}{\spec K}$ est
une courbe lisse et connexe
de $\P{K}{3}$, de degr\'e $d$ et genre g. Puisque
$g>G(d,n+1)$, alors, par d\'efinition de $G(d,n+1)$, il existe une surface
$X\!\subset\!\P{K}{3}$ de
degr\'e $n$ (\'eventuellement r\'eductible) contenant ${\ee C}_b$. Soit
$\xi\!\in\!\cohr{0}{\OO_{\P{K}{3}}(n)}$ le polyn\^ome homog\`ene
d\'efinissant
$X$. Posons $r\hs{-.5}=\hs{-1}-{\rm min}\{v(c),\ c$ est un coefficient de
$\xi
\}$ et $h=\xi{\pi^r}$, alors $h$ est un polyn\^ome non nul \`a coefficients
dans $A$, car
$v(c.\pi^r)=v(c)+r\geq 0$ pour tout coefficient $c$ de $\xi$, d\'efinissant
ainsi une surface
$\ee T$ de degr\'e $n$ dans $\P{A}{3}$. Nous allons voir que $\ee T$
r\'epond aux assertions (i) et
(ii).\par
\hs{5}(i) $\ $ Les deux fibres du morphisme induit ${\ee T}\to \ZZ$ \'etant
d\'efinies par des
polyn\^omes non nuls, car la fibre ${\ee T}_a$ (resp. ${\ee T}_b$)
au-dessus de
$a$ (resp. $b$) est d\'efinie par la classe de $h$ modulo l'id\'eal maximal
de $A$ (resp. \'egale
\`a $X$), alors ([M2], Example P, p. 299) ${\ee T}\to\ZZ$ est plat. Le
morphisme ${\ee C}\to\ZZ$
\'etant plat, on a $\ee C$ irr\'eductible et son point g\'en\'erique $\eta$
s'envoie sur $b$. Comme
${\ee C}_b$ est contenue dans ${\ee T}_b$, alors
$\eta\in{\ee T}$. On en d\'eduit l'inclusion ${\ee C}\subset {\ee T}$
sachant que ${\ee T}$ est un
sous-sch\'ema ferm\'e de $\P{A}{3}$.\smallskip
\hs{5}(ii)  D'apr\`es ce qui pr\'ec\`ede, on a $C\!=\!{\ee C}_a\!\subset\!
{\ee T}_a\cap S$ et
deg$({\ee T}_a)\hs{-1}=\hs{-1}n$, ainsi l'hypoths\`ese deg$(C)>s.n$
implique  n\'ecessairement
$\dim ({\ee T}_a\cap {S})\!=\!2$; comme $S$ est irr\'eductible, ceci
implique
$S\!\subset\!{\ee T}_a$. Soit $f$ (resp. $h_a$) le polyn\^ome d\'efinissant
$S$, (resp.$\ {\ee
S}_a$), alors $f$ divise $h_a$, et la surface $\Lambda $ cherch\'ee est
celle d\'efinie par le
quotient ${\dis{h_a}\over\dis f}$.\cqfd\vs{3}
 
\ind Nous arrivons maintenant au r\'esutat principal de cette section. (
Nous verrons son utilit\'e,
dans la prochaine section, pour construire des exemples de composantes non
r\'eduites du sch\'ema
de Hilbert $\Hdg$.). Rappelons que la normalis\'ee $\st_0$ de $S_0$ est
suppos\'ee lisse, donc ce
n'est autre que l'\'eclatement de $S_0$ le long de $L_0$, c'est aussi la
transform\'ee de
$S_0$ via l'\'eclatement de $\P{}{3}$ le long de $L_0$. \vs{2}
\sec{Proposition}{Soient $S_0$, $L_0$, $C_0$, $r$, $s$ comme dans la
proposition (2.8). On note $d$
et
$g$ respectivement le degr\'e et le genre de $C_0$. Soit $\ct_0\subset
\st_0$ la transform\'ee 
stricte de $C_0$. Si\par   (1) $r>2s(s-1)^2$\par (2) $d>s\sigma$ et
$g>G(d,\sigma+1)$ pour un
certain entier $\sigma$ tel que $s<\sigma<2s$.\par  (3)
$\dcohr{0}{\OO_{\st_0}(-\w
C_0)(\sigma)}=\dcohr{1}{\OO_{\st_0}(-\w C_0)(s)}=0$.\par
  Alors, pour toute d\'eformation plate $\ee C$ dans
$\P{}{3}$ de $C_0$ via $A=k[[t]]$, il existe une surface ${\ee
S}\subset\P{A}{3}$ de degr\'e 
$\sigma$ v\'erifiant les assertions (i) et (ii) de la proposition (2.9).}
\dem Soit ${\ee C}\subset\P{A}{3}$ une d\'eformation plate de $C_0$. Comme
$d>s\sigma$ et
$g>G(d,\sigma+1)$ alors d'apr\`es (2.10), il existe une surface ${\ee
T}\subset \P{A}{3}$ de
degr\'e
$\sigma$, plate sur $A$ et contenant $\ee C$ comme sous-sch\'ema ferm\'e;
de plus la fibre
sp\'eciale de ${\ee T}$ est
\'egale \`a $S_0\cup \Lambda_0$ o\`u $\Lambda_0$ est une surface de degr\'e
$\sigma-s$.\par
 Il nous suffit de montrer que la surface ${\ee T}$ ainsi d\'efinie est
\'egale \`a ${\bf
\Lambda}\cup \ee S$ o\`u $\ee S$ (resp. $\bf\Lambda$) est une d\'eformation
plate, param\'etr\'ee par
$A$, de $S_0$ (resp.
$\Lambda_0$). En effet, puisque ${\ee C}$ est irr\'eductible et ${\ee C}
\subset {\ee S}\cup
\bf\Lambda
$, on a ${\ee C}\subset {\ee S}$ ou
 ${\ee C}\subset\bf\Lambda $. Cette derni\`ere inclusion est exclue
puisque, d'une part
aucune composante irr\'eductible de ${\Lambda_0}$ n'est contenue dans $S_0$
(car $S_0$ est
irr\'eductible), et $S_0\not\subset\Lambda_0$ (car deg$\Lambda_0<s$) et
donc $S_0\cap \Lambda_0$ est
une courbe de degr\'e $s(\sigma-s)$. D'autre part, comme deg $C_0
>s\sigma$, on a
${\ee C}_a=C_0\not\subset\Lambda_0$ ($a$ d\'esigne le point ferm\'e de
$\spec A$). D'o\`u
${\ee C}\not\subset\bf\Lambda$, et donc ${\ee C}\subset {\ee S}$. Ainsi
${\ee S}$ v\'erifie la
conclusion (i) de (2.9); en utilisant (2.8) d\'eduit l'assertion (ii) et 
la proposition est d\'emontr\'ee.\par
\ind  Soit $h\in\cohr{0}{\OO_{\P{A}{3}}(\sigma)}$ l'\'equation de la
surface ${\ee T}$, on a $
h=\sum_{i\geq 0}h_it^i$ o\`u
$h_i\in\cohr{0}{\OO_{\P{k}{3}}(\sigma)}$ pour tout $i\geq 0$. Pour tout $
n\geq 0$, on pose
${\ee T}_n=\fibp{\ee T}{\spec A}{D_n}\subset \P{D_n}{3}$ et
$C_n=\fibp{\ee C}{\spec A}{D_n}\subset \P{D_n}{3}$.\par
 Pour tout $n\geq 0$, soit $h_{(n)}$ l'\'equation de ${\ee T}_n$ dans
$\P{D_n}{3}$; par
d\'efinition,
$h_{(n)}=h_0+h_1\eps{n}{}+\dots +\eps{n}{n}h_n$. \par Avec ces notations,
pour montrer l'existence
des surfaces $\bf\Lambda$ et $\ee S$ comme ci-dessus, il revient au m\^eme
de d\'emontrer le
\par
\sec{Lemme}{Il existe deux suites
$(f_{(n)})_n\in{\dis\prod_{n\geq 0}}\cohr{0}{{\OO}_{\P{D_n}{3}}(s)}$ 
\ et\ \ $(\lambda_{(n)})_n\in{\dis\prod_{n\geq
0}}\cohr{0}{{\OO}_{\P{D_n}{3}}(\sigma-s)}$ telles que
pour tout $n\geq 0$ on ait\par
\hs{10} (i)  $f_{(n)}\mid_{D_{n-1}} = f_{(n-1)}$\  ,\ 
$\lambda_{(n)}\mid_{D_{n-1}} =
\lambda_{(n-1)}$\par 
\hs{10} (ii) $h_{(n)}\ =f_{(n)}\lambda_{(n)}$\ .}
\ind En effet, les suites $(f_{(n)})_n$ et $(\lambda_{(n)})_n$ du lemme
sont, gr\^ace \`a la
propri\'et\'e $(i)$,  respectivement des \'el\'ements de
$\liminv_n\cohr{0}{{\OO}_{\P{D_n}{3}}(s)}$ et
${\liminv_{n}}\cohr{0}{{\OO}_{\P{D_n}{3}}(\sigma-s)}$. Par ailleurs,
puisque $A\simeq
\liminv_n\! D_n$, on\smallskip a
$\liminv_n\cohr{0}{\OO_{\P{D_n}{3}}(j)}\simeq
\cohr{0}{\OO_{\P{A}{3}}(j)}$ pour tout entier $j$. Soit
$f=\liminv_nf_{(n)}$ et $\lambda=\liminv_n\lambda_{(n)}$, alors
la\smallskip propri\'et\'e $(ii)$
implique
$h=f\lambda$. On pose $\ee S$ (resp. $\bf\Lambda$) la surface de $\P{A}{3}$
d\'efinie par
$f$ (resp. $\lambda$), et le lemme est d\'emontr\'e.\vs{3}  
{\it ( Preuve du lemme).} On fait une r\'ecurrence sur $n$. Pour
$n=0$, on a ${\ee T}_0= S_0\cup\Lambda_0$. Notons $f_{(0)}$ (resp.
$\lambda_{(0)}$) l'\'equation de
$S_0$ (resp. $\Lambda_0$). Posons $f_{(0)}=f_0$ et
$\lambda_{(0)}=\lambda_0$, on a alors
$h_{(0)}=f_{(0)}\lambda_{(0)}$.\par \ind Supposons le lemme vrai jusqu'\`a
l'ordre $n$. Il existe
alors $f_{(n)}\in\cohr{0}{{\OO}_{\P{D_n}{3}}(s)}$  et 
$\lambda_{(n)}\in\cohr{0}{{\OO}_{\P{D_n}{3}}(\sigma-s)}$ tels que
$h_{(n)}=f_{(n)}\lambda_{(n)}$.
 Soient $f_1,\dots ,f_n$ (resp. $\lambda_1,\dots,\lambda_n$ ) des
\'el\'ements de $\cohr{0}{{\OO}_{\P{k}{3}}(s)}$ (resp.
$\cohr{0}{{\OO}_{\P{k}{3}}(\sigma-s)}$ ) 
 tels que:
 $$f_{(n)}=f_0+\eps{n}{}f_1+\cdots +\eps{n}{n}f_n\ ,\ \
\lambda_{(n)}=\lambda_0+\eps{n}{}\lambda_1+\cdots +\eps{n}{n}\lambda_n\ .$$
Comme
$h_{(n+1)}\mid_{D_n}=h_{(n)}$ alors il existe
$G\in\cohr{0}{{\OO}_{\P{k}{3}}(\sigma)}$ tel que 
$$h_{(n+1)}=(f_0+\eps{n+1}{}f_1+\cdots +\eps{n+1}{n}f_n)\
(\lambda_0+\eps{n+1}{}\lambda_1+\cdots
+\eps{n+1}{n}\lambda_n)+
\eps{n+1}{n+1}G\ .$$ Il suffit de prouver que $G$ est dans l'id\'eal
$(f_0,\lambda_0)$. En effet, si
$G=af_0+b\lambda_0$ o\`u
$a\in\cohr{0}{\OO_{\P{k}{3}}(\sigma-s)}$ et
$b\in\cohr{0}{\OO_{\P{k}{3}}(s)}$, alors\smallskip
\hs{10}$h_{(n+1)}=(f_0+\eps{n+1}{}f_1+\cdots
+\eps{n+1}{n}f_n+\eps{n+1}{n+1}b)\ 
(\lambda_0+\eps{n+1}{}\lambda_1+\cdots
+\eps{n+1}{n}\lambda_n+\eps{n+1}{n+1}a)\ .$\par En
posant$f_{(n+1)}=f_0+\cdots +\eps{n+1}{n}f_n+\eps{n+1}{n+1}b\ {\rm  et}\
\lambda_{(n+1)}=\lambda_0+\cdots +\eps{n+1}{n}\lambda_n+\eps{n+1}{n+1}a\ ,$
le lemme est
d\'emontr\'e pour le rang $n+1$, et donc pour tout $n\geq0$.
\vs{1}
\ind Pour prouver que $G$ est dans $(f_0,\lambda_0)$, nous allons
proc\'eder en plusieurs pas. \par
\ind Soit $T$ le sous-sch\'ema de $\P{k}{3}$ d\'efini par l'id\'eal
$(f_0,\lambda_0)$. Soit
$\varphi$ le morphisme d\'efini dans (2.1). On note $\w T=\varphi^{-1}T$ et
$\w G=\varphi^*G$, notons
ici que $\w G\in\cohr{0}{\OO_{\st_0}(\sigma)}$ et que
$\OO_{\st_0}(\w T)\simeq \OO_{\st_0}(\sigma-s)$. Le but de ce qui suit est
de montrer que
$G\mid_T = 0$, puique dans le
$2$\`eme pas ci-dessous on montrera que ceci suffit pour que $G$ soit dans
$(f_0,\lambda_0)$.
L'id\'ee consiste
\`a montrer que $\w G\mid_{\w T}=0$ (cf. pas 4), ensuite, et c'est
l'op\'eration la plus difficile,
passer de
$\w G\mid_{\w T}=0$ \`a $G\mid_T = 0$. Pour cela on verra qu'en fait $\w
G\mid_{\w T}=0$ implique
que $G$ s'annule sur $T$ sauf peut-\^etre aux points-pinces, et en
utilisant la preuve du $2$\`eme
pas, on montre que les points-pinces annulent $G$ car sinon ce sont des
points immerg\'es de $T'$,
le c\^one affine sur $T$, mais ceci est exclu en vertu du premier pas, donc
$G\mid_T=0$.
\vs{2}
\ind $Pas\ 1.$ Soit $T'$ le c\^one affine sur $T$, alors $T'$ n'a pas de
points immerg\'es. En
effet, puisque $f_0$ est irr\'eductible et
deg$\lambda_0$$=\sigma-s$$<s=$deg$f_0$, les polyn\^omes
$f_0$ et
$\lambda_0$ n'ont pas de composantes communes. Donc codim $T'$ = codim
$T=2$, ce qui implique que
$T'$ est l'intersection compl\`ete de
$V(f_0)$ et $V(\lambda_0)$, et d'apr\`es ([Ar], Th\'eor\`eme 5.1), il est
ainsi Cohen-Macaulay, et
donc d\'epourvu de points immerg\'es. \par 
\ind $Pas\ 2.$ Si $G\mid_T = 0$ alors $G\in(f_0,\lambda_0)$. En effet,
$G\mid_T = 0$ implique $
G\mid_{T'} =0$ sauf peut \^etre en $O$, l'origine de
$\A{k}{4}$. Si $G\not =0$ dans $\OO_{T',O}$ , alors $G$ est une section de
$\OO_{T'}$ support\'ee en $O$, et donc $O$ est un point immerg\'e de $T'$,
ce qui est exclu par le 
premier pas. Donc $G\mid_{T'} =0$, et comme
$\coh{0}{T'}{\OO_{T'}} =\coh{0}{\A{}{4}}{\OO_{\A{}{4}}}/_{\displaystyle
(f_0,\lambda_0)}$, alors
$G\in (f_0,\lambda_0)$.\par
\ind $Pas\ 3.$ Montrons $\w G\mid_{\w T\cap {\w {C_0}}} = 0$. Pour tout
$m\leq n$, notons $S_m$
(resp. $\Lambda_m$) la surface de $\P{D_m}{3}$ d\'efinie par $f_{(m)}$
(resp. $\lambda_{(m)}$). Par
hypoth\`ese de r\'ecurrence, on a
${\ee T}_n=S_n\cup\Lambda_n$, et on montre comme dans la preuve de (2.11)
que $C_n\subset S_n$. Donc
$(f_0+\cdots +\eps{n+1}{n} f_n)\mid_{C_n} =f_{(n)}\mid_{C_n}=0$; il existe
alors
$\beta\in\cohr{0}{\OO_{C_0}(s)}$ tel que $(f_0+\cdots +\eps{n+1}{n}
f_n)\mid_{C_{n+1}}=\eps{n+1}{n+1}\beta$. Par ailleurs on a $h\mid_{C_{n+1}}
=h_{(n+1)}
\mid_{C_{n+1}}=0$, car ${\ee C}\subset {\ee T}$, donc\smallskip
\hs{30} $\eps{n+1}{n+1}[ \beta (\lambda_0+\cdots +\eps{n+1}{n}\lambda_n) +
G]\mid_{C_{n+1}} =0\
$,\smallskip i.e., $ (\lambda_0\beta\ +\ G)\mid_{C_0}=0$. Ceci entra\^ine
$\ \w G\mid_{\w C_0}\in (\w \lambda_0\mid_{\w C_0})$, or $\w T\cap\w C_0$
est le sous-sch\'ema de
$\w C_0$ d\'efini par $\w\lambda_0\mid_{\w C_0}$, donc $\w G\mid_{\w T\cap
{\w {C_0}}}\ =0$.  Ceci
sera essentiel dans le pas suivant.\vs{2} 
\ind $Pas\ 4.$ On montre $\w G\mid_{\w T}=0$. Consid\'erons la suite
exacte\smallskip
\hs{40}$0\to\OO_{\st_0}(-\w C_0)(\sigma)\to\OO_{\st_0}(\sigma)\to\OO_{\w
C_0}(\sigma)\to
0$\smallskip que l'on tensorise par
$\OO_{\w T}$\ , on obtient alors la suite exacte\smallskip 
\hs{32}$0\to \HH\to\OO_{\w T}(-\w C_0\cap\w T)(\sigma )\to\OO_{\w
T}(\sigma)\to\OO_{\w T\cap\w
C_0}(\sigma)\to 0$\smallskip o\`u
$\HH =\Torr{\OO_{\st_0}}{1}{\OO_{\w C_0}}{\OO_{\w T}}$. Soit $\KK$ le noyau
de $\OO_{\w
T}(\sigma)\to\OO_{\w T\cap\w C_0}(\sigma)$, on a alors la suite exacte
\smallskip
\hs{40}$0\to\HH\to\OO_{\w T}(-\w C_0\cap\w T)(\sigma )\to \KK\to
0$\smallskip Comme $\HH$ est
support\'e par $\w T\cap\w C_0$ qui est un ensemble fini, alors
$H^1(\HH)=0$, d'o\`u la suite
exacte\smallskip 
\hs{35}$0\to H^0(\HH)\to H^0(\OO_{\w T}(-\w C_0\cap\w T)(\sigma ))\to
H^0(\KK)\to 0$\smallskip
\ind Consid\'erons maintenant, en rappelant $\OO_{\st_0}(\w T)\simeq
\OO_{\st_0}(\sigma-s)$, le
diagramme commutatif suivant: 
$$\diagram{&&&&0\cr &&&&\vfl{}{}\cr   0&\fh&\OO_{\st_0}(-\w
C_0)(s)&\fh&\OO_{\st_0}(-\w
C_0)(\sigma)&\fh&\OO_{\w T}(-\w C_0\cap\w T)(\sigma)&\fh&0\cr 
&&\vfl{}{}&&\vfl{}{}&&\vfl{}{}\cr
0&\fh&\OO_{\st_0}(s)&\fh&\OO_{\st_0}(\sigma)&\fh&\OO_{\w
T}(\sigma)&\fh&0\cr
&&&&\vfl{}{}&&\vfl{}{}\cr &&&&\OO_{\w C_0}(\sigma)&\fh&\OO_{\w T\cap\w
C_0}(\sigma\cr
&&&&\vfl{}{}&&\vfl{}{}\cr &&&&0&&0\cr}$$
 D'apr\`es le pas 3, on a $\w G\mid_{\w T\cap {\w
{C_0}}}=0$, alors
$\w G\in \cohr{0}{\OO_{\st_0}(\sigma)}$ s'envoie sur $0$ dans
$\cohr{0}{\OO_{\w T\cap\w
C_0}(\sigma)}$, ainsi $\w G\!\mid_{\w T} \in H^0(\KK)$. Par ailleurs, on a
$\cohr{0}{\OO_{\w T}(-\w C_0\cap\w T)(\sigma)}=0$ car, par hypoth\`eses,
$\dcohr{0}{\OO_{\st_0}(-\w C_0)(\sigma)}=\dcohr{1}{\OO_{\st_0}(-\w
C_0)(s)}=0$, donc $H^0(\KK)=0$, et
donc $\w G\mid_{\w T}=0$.  \vs{2}  
\ind Notons $\w \lambda_0\!=\!\varphi^*(\lambda_0)$, c'est une section de
$\OO_{\st_0}(\sigma-s)$ dont $\w T$ est le sch\'ema des z\'eros. D'apr\`es
le pas pr\'ec\'edent, on
a $\w G\mid_{\w T} = 0$, donc il existe une section $\w \alpha\in
\cohr{0}{\OO_{\st_0}(s)}$ telle
que  ${\w G}= {\w\lambda_0} {\w\alpha}$. Notons $\PP$ l'ensemble des
points-pinces. Nous allons
montrer \`a laide des pas qui suivent qu'il existe une section
$\alpha'\!\in\!\cohr{0}{\OO_{S'_0}(s)}$ telle que $\w
\alpha =\varphi^*\alpha'$ sur $\w S'_0$, o\`u $S'_0=S_0\setminus\PP$ et $\w
S'_0=
\varphi^{-1}(S'_0)$. \vs{2}
\ind $Pas\ 5.$ On se place dans le voisinage \'etale $(V_n,e_n)$ dans
$\P{D_n}{3}$ des points de
$C_0\cap L_0$, d\'efini dans (2.5). Puisqu'on est dans les hypoth\`eses de
la proposition (2.6), il
existe des param\`etres uniformisants $u_n$ et $v_n$ de $V_n$ tels que
$e_n^*f_{(n)}=u_nv_n$. ( Pour
\'eviter d'allourdir les notations, et tant que le contexte est clair, on omettra les $e_n^*$). 
Posons
$F=f_0+\cdots +\varepsilon_{n+1}^n f_n$, et soient $u$ et $v$ des sections
de
$\cohr{0}{\OO_{V_{n+1}}(s)}$ telles que $u\mid_{V_n}=u_n$ et
$v\mid_{V_n}=v_n$. Il existe
$\beta\in\cohr{0}{\OO_{V_{n+1}}(s)}$ tel  que 
 $F = uv +\eps{n+1}{n+1}\beta$.\par
 \ind Notons $L^*_0$ l'image inverse de $L_0$ dans $V_{n+1}$ par le
morphisme $e_{n+1}$, alors
$\beta\mid_{L^*_0}$ ne d\'epend pas des rel\`evements choisis
 $u$ et $v$. En effet, soient $u'$ et $v'$ deux autre sections telles que
$u'\mid_{V_n}=u_n$ et
$v'\mid_{V_n}=v_n$ , il existe alors $a$, $b$, $\beta'$ dans
$H^0(\OO_{V_{n+1}}(s))$ tels que
$$u' =u +\eps{n+1}{n+1}a\ , \ v' =v + \eps{n+1}{n+1}b \ ,
\ F = u'v' +\eps{n+1}{n+1}\beta'\ .$$ De l'\'egalit\'e $u'v'+
\eps{n+1}{n+1}\beta'= 
uv+\eps{n+1}{n+1}\beta$ on tire $\beta=\beta'+av+bu\mod{\eps{n+1}{}}$, et
donc 
$\beta\mid_{L^*_0} =\beta'\mid_{L^*_0}$ puisque $L^*_0$ est d\'efini dans
$V_{n+1}$ par l'id\'eal $(u,v,\eps{n+1}{})$. Ceci permet de voir que
$\beta\mid_{L^*_0}$ descend en
une section rationnelle de $\OO_{L_0}(s)$, qu'on note $\beta_{L_0}$, qui
est r\'eguli\`ere sauf
peut-\^etre aux points-pinces.\vs{2}
\ind $Pas\ 6.$ On pose $\w \beta= \varphi^*(\beta_{L_0})$, alors \#(p\^oles
de
$\w\beta) \leq 2[s(s-1)^2-\deg \OO_{L_0}(s)]$. En effet, par une variante
du th\'eor\`eme
de Bertini, pour une section g\'en\'erale $f_{n+1}$ de $\OO_{\P{}{3}}(s)$,
les z\'eros de
$f_{n+1}\mid_{L_0}+\beta_{L_0}$ sont tous simples et ses p\^oles sont ceux
de $\beta_{L_0}$. Par
ailleurs,  en tout point de $L_0$ o\`u $f_{n+1}\mid_{L_0}+\beta_{L_0}$
s'annule, le sch\'ema
$${\bf Pol}( q:\w\P{D_{n+1}}{3}\hfl{proj.g\acute
en.}{}\P{D_{n+1}}{1},F+\eps{n+1}{n+1}f_{n+1})$$  a
un point immerg\'e (cf. preuve du sous-lemme (2.7.1)). Donc \#$\{$z\'eros
de
$(f_{n+1}\mid_{L_0}+\beta_{L_0})\} \leq  s(s-1)^2$. Or \#$\{$p\^oles de
$\w\beta\}= 2$\#$\{$p\^oles
de $\beta_{L_0}\}$, car $\deg\varphi=2$. Donc
\#$\{$p\^oles de $\w\beta\}= 2$\#$\{$z\'eros de
$(f_{n+1}\mid_{L_0}+\beta_{L_0})\}-2\deg\OO_{L_0}(s)\leq 2[s(s-1)^2-\deg
\OO_{L_0}(s)]$.\vs{1} 

\ind $Pas\ 7.$ On montre que pour tout $ P\in \w C_0\cap\w L_0$, on a $\w
\alpha (P)\ =- \w \beta
(P)$. Soit $P\in\w C_0\cap\w L_0$, notons $Q$ son image par $\varphi$.
Comme dans le pas 5, on se
place dans un voisinage \'etale de $Q$. On a $F= uv +
\eps{n+1}{n+1}\beta$, et $ h_{n+1}=F\Lambda+\eps{n+1}{n+1}G$, o\`u $\Lambda
=\lambda_0+\cdots
+\eps{n+1}{n}\lambda_n$.\par
\ind Puisque $\eps{n+1}{n+1}\mid_{C_n}=0$, $h_{n+1}\mid_{C_n}=0$ et
$\Lambda\mid_{C_0}\not =0$, on a
$F\mid_{C_n}=0$, et donc $ uv\mid_{C_n}=0$. On suppose, par exemple, que la
courbe 
$C_0$ est contenu dans la nappe $(u=0)$. Puisque $C_0$ est transverse \`a
${L_0}$, 
$v\mid_{C_0}$ n'est pas diviseur de z\'ero, donc $v\mid_{C_n}$ ne l'est pas
non plus. D'o\`u
$u\mid_{C_n}=0$. (Dans l'autre cas $v\mid_{C_n}=0$). Ainsi il existe
$\gamma\in\cohr{0}{\OO_{C_0}(s)}$ tel que
$u\mid_{C_{n+1}}=\eps{n+1}{n+1}\gamma$. Par
ailleurs,\smallskip puisque $h_{n+1}\mid_{C_{n+1}} =F\mid_{C_{n+1}}
\Lambda\mid_{C_{n+1}}+ \eps{n+1}{n+1}(G\mid_{C_0})=0$, on a\vs{2}
\hs{34} $[u\mid_{C_{n+1}}v\mid_{C_{n+1}} + \eps{n+1}{n+1}(\beta\mid_{C_0})]
\Lambda\mid_{C_{n+1}}+
\eps{n+1}{n+1}(G\mid_{C_0})=0.$\par Donc\par \hs{40}$
\eps{n+1}{n+1}[G\mid_{C_0}\
+\Lambda\mid_{C_0}(\gamma\ v\mid_{C_0}\ +\ \beta\mid_{C_0})]=0$.\par
 Ce qui implique $$G\mid_{C_0} + (\gamma\mid_{C_0}\ v_0\mid_{C_0}
+\beta\mid_{C_0})\
\lambda_0\mid_{C_0} =0.$$ D'o\`u l'\'egalit\'e
$ {{G\mid_{C_0}}\over\lambda_0\mid_{C_0}} =- (\gamma v_0
+\beta)\mid_{C_0}$, car
$\lambda_0\mid_{C_0}\not =0$. Mais $\w\alpha\mid_{\w C_0} = {\w G\mid_{\w
C_0}\over\w
\lambda_0\mid_{\w C_0}}$ et $\w v_0(P)=0$, donc $ \w \alpha (P) =-\w \beta
(P)\ $.\vs{2}
\ind $Pas\ 8.$ Notons $\PP$ l'ensemble des points-pinces. Nous allons
montrer qu'il existe une
section $\alpha'\!\in\!\cohr{0}{\OO_{S'_0}(s)}$ telle que $\w
\alpha =\varphi^*\alpha'$ sur $\w S'_0$, o\`u $S'_0=S_0\setminus\PP$ et $\w
S'_0=
\varphi^{-1}(S'_0)$.\par
\ind Comme $\w \alpha\mid_{\w L}$ est une section (r\'eguli\`ere) de
$\OO_{\w L_0}(s)$ (qui est un fibr\'e de degr\'e $2\deg\OO_{L_0}(s) $ ) et
$\w \beta+ \w
\alpha\mid_{\w L_0}$ est une section rationnelle de $\OO_{\w L_0}(s)$ dont
les p\^oles sont ceux de
$\w \beta$, alors d'apr\`es le pas 6, on a \#$\{$ z\'eros de $(\w \beta+\
\w \alpha\mid_{\w
L_0})\}\leq 2s(s-1)^2$. Or par hypoth\`ese $\#(\w C_0\cap \w L_0) = \#
(C_0\cap L_0)>2s(s-1)^2$,
donc d'apr\`es le pas 7, $ \w\beta+\w \alpha\mid_{\w L_0}=0$. Ainsi
$\w\beta$ est r\'eguli\`ere (en
particulier
$\beta_{L_0}$ est d\'efini sur $L_0$) et $\w
\alpha\mid_{\w L_0}=-\varphi^*\beta_{L_0}$ en dehors des points de
ramification. Nous allons voir
que $\w\alpha\mid_{\w L_0}$ descend en une section sur $S'_0$.\par
\ind Voyons d'abord ceci g\'eom\'etriquement; $\w S'_0$ (qui est la
normalis\'ee de $S'_0$) est
localement la r\'eunion de deux nappes, et en recollant ces nappes (on
recolle en fait les paires de
l'involution induite sur $\w L'_0=\w L_0\setminus \varphi^{-1}\PP$), on
obtient la surface $S'_0$.
Pour avoir  une section sur $S'_0$, il faut avoir une section sur $\w S'_0$
qui a la m\^eme valeur
sur les deux nappes. Or en dehors des points de ramification on a $\w
\alpha\mid_{\w
L_0}=-\varphi^*\beta_{L_0}$, donc $\w\alpha\mid_{\w L_0}$ est fix\'ee par
l'involution induite sur
$\lt_0$. Donc $\w\alpha\mid_{\w L_0}$ descend en une section sur
$L'_0$, et ceci implique que $\w\alpha$ descend en une section $\alpha'$
sur $S'_0$.\par
\ind Une autre fa{\c c}on de voir cela consiste \`a consid\'erer le
diagramme commutatif (de
normalisation) suivant:
$$\diagram{0&\hfl{}{}&\OO_{S'_0}(s)&\hfl{}{}&\varphi_*\OO_{\w
S'_0}(s)&\hfl{}{}&\FF&\hfl{}{}&0\cr
&&\vfl{}{}&&\vfl{}{}&&\parallel\cr
0&\hfl{}{}&\OO_{L'_0}(s)&\hfl{}{}&\varphi_*\OO_{\w
L'_0}(s)&\hfl{}{}&\FF&\hfl{}{}&0\cr}$$  o\`u $\FF$ est un faisceau de rang
1 localement libre sur
$L'_0$. On obtient le diagramme commutatif
suivant:$$\diagram{0&\hfl{}{}&H^0(\OO_{S'_0}(s))&\hfl{\varphi^*}{}&H^0(\OO_{
\w
S'_0}(s))&\hfl{a}{}&H^0(\FF)&\hfl{}{}&0\cr
&&\vfl{}{}&&\vfl{}{}&&\parallel\cr
0&\hfl{}{}&H^0(\OO_{L'_0}(s))&\hfl{\varphi^*}{}&H^0(\OO_{\w
L'_0}(s))&\hfl{b}{}&H^0(\FF)&\hfl{}{}&0\cr}$$  On obtient $\alpha'$ en
observant $\ 0=b\circ
\varphi^*(\beta_{L_0})=b(\w\alpha\mid_{\w L'_0})=a(\w\alpha)$ \par
\ind {\it Conclusion. } On a $G\mid_{ S'_0}=\alpha'
\lambda_0\mid_{ S'_0}$, donc $G\mid_{T^*}=0$ o\`u $T^*= T\setminus \PP$. Or
$T$ est une intersection
compl\`ete de dimension  $1$, donc $T$ n'a pas de points immerg\'es, et un
argument analogue \`a
celui du pas 2 on d\'eduit $G\mid_T =0$, ce qu'on cherchait.\cqfd 
 \parag{Application au sch\'ema de Hilbert $\Hdg$.}
\ind Dans cette section, le r\'esutat principal que nous allons prouver est
le\par
\sec{Th\'eor\`eme}{ Soient $d$ et  $g$ deux entiers satisfaisant:\smallskip
$(A)$\hs{30}\ \ $G(d,8)<g\leq {\dis {73\over 2}}(d-74) $\ \  si\ \ $146\geq
d\geq 95$\smallskip
$(B)$\hs{30} \ \  $G(d,8)<g\leq {\dis{(d-1)^2\over 8}}$\ \  \ si\ \
$d\geq147$\smallskip  Alors il
existe une composante irr\'eductible non r\'eduite du sch\'ema de Hilbert
$\Hdg$ dont l'\'el\'ement
g\'en\'eral est une courbe trac\'ee sur une surface quartique \`a droite
double. }  
\sec{Notations et g\'en\'eralit\'es }{} 
{\bf 3.2.1.  }On note $H(4)$ l'espace projectif
$\P{}{}\cohr{0}{\OO_{\P{}{3}}(4)}$ des surfaces quartiques dans $\P{}{3}$.
Soit $D(d,g;4)\subset
\fibp{\Hdg }{}{H(4)}$ le sch\'ema de Hilbert des drapeaux courbes-surfaces.
Pour toute
droite
$L$ de
$\P{}{3}$ on pose
 $\QQ_L=\P{}{}\cohr{0}{\JJ^2_{L/\P{}{3}}(4)}$; ensemblistement $\QQ_L$ est
l'espace des quartiques
contenant $L$ comme droite double.\par
 Soit $q\colon D(1,0,4)\to H_{1,0}$ la projection naturelle, associant \`a
chaque paire $(L,S)$ la
droite $L$. La fibre de $q$ au-dessus de chaque $L\in H_{1,0}$ \'etant
\'egale \`a
$\P{}{}\cohr{0}{\JJ_{L/\P{}{3}}(4)}$, on peut identifier
$D(1,0,4)$ \`a $\P{}{}(p_*\JJ_{\bf I}(4))$ o\`u ${\bf I}
\subset\P{}{3}\times H_{1,0}$ est la
vari\'et\'e d'incidence et $p$ est la projection $\P{}{3}\times H_{1,0}\to
H_{1,0}$. Posons
$\QQ':= \P{}{}(p_*\JJ^2_{\bf I}(4))$. Ensemblistement, $\QQ'$ est l'espace des paires
$(L,S)$ o\`u $L$ est une droite et $S$ est une surface quartique de
$\P{}{3}$, contenant $L$ comme
droite double. Soit $\QQ$ l'ouvert de $\QQ'$ form\'e par les couples
$(L,S)$ tels que $S$ est une
surface int\`egre, i.e., irr\'eductible et r\'eduite.  
\par
{\bf 3.2.2.  } Soit $\WW_{d,g}$ le sous-sch\'ema de $H_{1,0}\times
D(d,g,4)$ image r\'eciproque de
$\QQ$ par le morphisme de projection $H_{1,0}\times D(d,g,4)\to
H_{1,0}\times H(4)$. Ensemblistement,
$\WW_{d,g}$ est form\'e des triplets $(L ,C,S)$ o\`u $C$ est une courbe
lisse et connexe de degr\'e
$d$ et genre $g$ trac\'ee sur une surface quartique int\`egre contenant $L$
comme droite double. On
pose $\rho:\WW_{d,g}\to\Hdg $ la projection induite sur $\WW_{d,g}$. On
note $\VV$ l'image
sch\'ematique de $\WW_{d,g}$ par $\rho$, c'est la fermeture sch\'ematique
de l'ensemble des courbes
lisses et connexes de degr\'e
$d$ et de genre $g$ dans $\P{}{3}$, trac\'ees sur des surfaces quartiques
int\`egres \`a droite
double.\par
\ind Afin de donner un sens \`a la d\'efinition de $\WW_{d,g}$, on suppose
dans toute cette section
que $0\leq g\leq (d-1)^2/8$, car ceci assure  d'apr\`es le th\'eor\`eme de
Gruson-Peskine ([GP2],
th\'eor\`eme 1.1) que $\WW_{d,g}$ est non vide et domine $\QQ$. \vs{3}
{\bf 3.2.3.  } Dans toute cette section on d\'esignera par
$S_0$ une surface quartique int\`egre dans $\P{}{3}$, ayant une droite
double $L_0$. Soit
$\pi\colon\w\P{}{3}\to\P{}{3}$ l'\'eclatement de
$\P{}{3}$ le long de $L_0$. Soit $\st_0$ la transform\'ee stricte de $S_0$
par $\pi$. Notons
$\varphi$ la restriction $\pi\hs{-1}\mid_{\st_0}\ $. On choisit $(L_0,S_0)$
suffisemment g\'en\'erale
dans
$\QQ$, dans ce cas la surface $\st_0$ (qui est aussi la normalis\'ee de
$S_0$) est lisse (cf.
[Az], I-1.3).  Ainsi en appliquant ([S], 1.1), $\st_0$ est isomorphe \`a
l'\'eclatement de $\P{}{2}$
en 9 points situ\'es sur une cubique lisse. On a alors Pic $\st_0\simeq
{\bb Z}^{10}$, et on peut
prendre comme base orthogonale, la famille $\{\Delta, -E_1,\dots,-E_9\}$
o\`u $\Delta$ est l'image
inverse d'une droite g\'en\'erale de $\P{}{2}$ et o\`u les $E_i$ sont les
classes des droites
exceptionnelles. Une telle base est dite {\it exceptionnelle}. Le morphisme
$\varphi\colon\st\to\P{}{3}$ est donn\'e par le syst\`eme lin\'eaire
$H=\mid(4,2,1^8)\mid$, la
classe $K_{\st_0}$ du fibr\'e canonique de
$\st_0$ \'egale \`a $(-3,-1^9)$, et la classe de $\lt_0=\varphi^{-1}L_0$
est \'egale \`a
$-K_{\st_0}$. Le morphisme $\varphi$ induit un morphisme $\lt_0\fh L_0$ de
degr\'e $2$ ramifi\'e en
$4$ points que l'on notera $R_1,\dots , R_4$. ( Pour les d\'etails sur
cette construction, on renvoit
\`a  [GP2] ). Signalons enfin que les hyoth\`eses faites sur $S_0$
impliquent, en vertu de
([R], 4.1), que c'est une surface \`a singularit\'es ordinaires.\vs{5}
\ind Dans le th\'eor\`eme suivant nous donnons des conditions suffisantes
pour qu'une composante
irr\'eductible de $\Hdg$, contenant une courbe $C_0$ de degr\'e $d$ et
genre
$g$ trac\'ee sur $S_0$, soit non r\'eduite.\par
\sec{Th\'eor\`eme}{Soit $S_0$ une surface quartique \`a droite double $L_0$
comme dans (3.2.3). Soit
$C_0\subset S_0$ une courbe lisse et connexe de degr\'e $d$ et genre $g$
distincte de $L_0$, on
note $\ct_0$ sa transform\'ee stricte dans $\st_0$. Soit
$V$ une composante irr\'euctible $\Hdg$ contenant
$C_0$. On suppose qu'il existe une composante irr\'eductible
$\WW$ de ${\WW}_{d,g}$ qui domine $V$. Si
$h^0(\OO_{\st_0}(4H-\ct_0-\lt_0))=0$ et $ h^1{(\NN_{\ct_0/\st_0})}=0 $,
alors $V$ est
g\'en\'eriquement non r\'eduite.}
\ind L'id\'ee de la preuve consiste \`a montrer l'in\'egalit\'e dim
${T_{C}H_{d,g}}>\dim V$
o\`u $C$ est l'\'el\'ement g\'en\'erique de $V$. Par semi-continuit\'e, il
suffit de montrer cette
in\'egalit\'e en $C_0$, et comme dim $V\leq$ dim
$\WW\leq $ dim $T_{(L_0,C_0,S_0)}\WW$, il suffit alors de voir, en notant
$t$ la dimension de $T$, 
que $t_{(L_0,C_0,S_0)}\WW<t_{C_0}(H_{d,g})$. Pour ce faire nous aurons
besoin des r\'esultats
interm\'ediaires suivants:  \smallskip
\sec{Proposition}{En utilisant les notations et les hypoth\`eses de (3.2.3), on consid\`ere le
faisceau de $\OO_{\st_0}$-modules $\MM$ conoyau du morphisme tangent
$T_\varphi\colon\tst\to\fitp$.
Soient
$\NN={\MM}^{\vee\vee}$ le double dual de $\MM$ et $\delta\colon\MM\to\NN$
le morphisme canonique,
alors:\par 
\ind 1) $\MM$ est sans torsion, donc en particulier $\delta$ est
injectif.\par
\vskip -2pt \ind 2) Le conoyau de $\delta$ est isomorphe \`a $\cq$, o\`u
$k(R_i)$ d\'esigne le
faisceau gratte-ciel au-dessus de $R_i$.\par
\ind 3) ll existe un isomorphisme $\NN\simeq\nstpt(\lt_0)$, o\`u $\nstpt$
est le fibr\'e normal de
$\st_0$ dans $\w\P{}{}$.}
\dem 1) Le noyau de $T_\varphi$ ayant pour support un sous-ensemble de $\{
R_1,\dots, R_4\}$ (cf.
[H3], 3.1) est donc de torsion, mais comme c'est un sous-faisceau d'un
fibr\'e vectoriel, il est
alors nul. Ainsi, $\MM$ est un faisceau coh\'erent, localement libre en
dehors de
$\{ R_1,\dots, R_4\}$ (qui est un ferm\'e de $\st_0$ de codimension
sup\'erieure \`a $2$), donc $\MM$
est  sans torsion.\par 2) Le faisceau $\NN$ \'etant r\'eflexif (cf. [H4],
1.2) de rang \'egal \`a
$1$ (celui de $\MM$) est donc inversible (cf. [H4], 1.9). Soit
$\eta\colon\NN\to\OO_{\st_0}$ un
isomorphisme, et consid\'erons la surjection canonique
$\psi\colon\fitp\to\MM$. Notons
$\upsilon\colon\OO_{\st_0}\to\cq$ la somme directe des morphismes
d'\'evaluation. Alors il suffit de
montrer que la suite 
$$0\to\tst\hfl{T_\varphi}{}\fitp\hfl{\theta}{}\NN\hfl{\tau}{}\cq\to0\leqno
(3.4.1)$$ est exacte, avec
$\tau=\upsilon{\circ}\eta$ et $\theta=\delta{\circ} \psi$.\par 
\ind Pla{\c c}ons-nous d'abord au voisinage d'un point de $\st_0$ qui n'est
pas de ramification,
alors
$\tau=0$ et $\delta$ est un isomorphisme; l'exactitude de la suite
$(3.4.1)$ est d\'eduite de celle
de la suite $$0\to\tst\hfl{T_\varphi}{}\fitp\hfl{\psi}{}\MM\to0.$$
\ind  Pla{\c c}ons-nous maintenant autour d'un point $R\in\{R_1,\dots ,
R_4\}$. Puisque $S_0$ est
\`a singularit\'es ordinaires, on sait qu'il existe des coordonn\'ees
locales $(s,t)$ de
$\st_0$ au voisinage de
$R$ telles que
$\varphi(s,t)=(s,st,t^2)$; ainsi la suite $(3.4.1)$ s'\'ecrit localement:
$$0\to k[s,t]^2\hfl{{}^tM}{}k[s,t]^3\hfl{Y}{} k[s,t]\hfl{r}{}k\to 0\leqno
(3.4.2)$$ \par o\`u
$M=\bigl[{1\atop 0}{t\atop s}{0\atop 2t}\bigr]$, $Y=[2t^2, -2t, s]$ et o\`u
$r$ est la r\'eduction
modulo l'id\'eal $(s,t)$ (pour trouver $Y$, on v\'erifie que c'est l'unique
matrice, \`a une
multiplication pr\`es par une fonction, telle $Y.{}^tM=0$.). On laisse au
lecteur le soin de
s'assurer que la suite $(3.4.2)$ est bien exacte.\par  3) Soit $\EE$ le
conoyau du morphisme tangent
$T_{\pi}\colon\tpt\hookrightarrow\pitp$, c'est un faisceau support\'e sur
la surface quadrique
$E=\pi^{-1}L_0$; on montre en fait que
$\EE\hs{-1}\simeq\hs{-.5}\OO_E(2H-E)$ (cf. [Az], II-3.2). Soient $T_{\mu}$
le morphisme tangent \`a
l'immersion ferm\'ee $\mu\colon{\st_0}\hook\pt$ et $\alpha$ la restriction
$T_\pi\hs{-1}\mid_{\st_0}$; alors il est facile de voir qu'il existe un
morphisme
$\xi\colon\nstpt\to\MM$ tel que le diagramme suivant soit commutatif

$$\diagram{ &&&&0&&0\cr &&&&\vfl{}{}&&\vfl{}{}\cr
0&\hfl{}{}&\tst&\hfl{T_{\mu}}{}&\tptst&\hfl{\gamma}{}&\nstpt&\hfl{}{}&0\cr
&&\bigver&&\vfl{\alpha}{}&&\vfl{\xi}{}\cr 
0&\hfl{}{}&\tst&\hfl{T_\varphi}{}&\fitp&\hfl{\psi}{}&\MM&\hfl{}{}&0\cr  
&&&&\vfl{\beta}{}&&\vfl{}{}\cr 
 &&&&\EE\mid_{\lt_0}&\simeq& Coker\ \xi\cr    &&&&\vfl{}{}&&\vfl{}{}\cr    
&&&&0&&0\cr}\leqno
(3.4.3)$$ o\`u $\beta$ et $\gamma$ sont les projections naturelles. Soit
$\zeta\colon\nstpt\to\NN$ le
compositum
$\delta{\circ}\xi$, et consid\'erons la projection
$\lambda\colon\NN\to\NN\mid_{\lt_0}$; comme $Ker\
(\lambda)\simeq\NN(-\lt_0)$, alors pour prouver $(3)$, il suffit de montrer
que la suite suivante
est exacte
$$0\to\nstpt\hfl{\zeta}{}\NN\hfl{\lambda}{}\NN\mid_{\lt_0}\to0\
.\leqno(3.4.4)$$
\ind Pour ce faire, on proc\`ede localement comme dans (2). C'est clair au
voisinage d'un point de
${\st_0}\setminus \lt_0$ (car dans ce cas on a
$\nstpt\simeq\MM\simeq\NN$ et $\NN\simeq\NN(-\lt_0)$).\par Si l'on se place
autour d'un point de
ramification, alors en utilisant les notations de (2), on voit
que\break\vs{-4.5} ${}^t[{1\atop
0}{0\atop 0}{0\atop 2t}]$,
$[0,-2t,1]$ et $\biggl[{{\sc1\atop \sc t }\atop\sc 0}{{\sc0\atop \sc s
}\atop\sc 0}{{\sc0\atop\sc 0
}\atop \sc1}\biggr]$ sont, respectivement, les matrices des morphismes
$T_{\mu}$, $\gamma$ et
$\alpha$ d\'efinis\vs{.5} dans le diagramme$(3.4.3)$. La commutativit\'e de
ce dernier implique
$\zeta\gamma=\delta\xi=\psi\alpha$; puisque la matrice de $\psi\alpha$ est
\'egale \`a $[0,-2st,s]$,
alors $\zeta$ a pour matrice $[s]$. Il reste
\`a observer que $\NN\mid_{\lt_0}$ est isomorphe \`a $k[t]$ et que
$\lambda$ est la r\'eduction
modulo
$(s)$.\par Pla{\c c}ons-nous maintenant en un point $R$ de $\lt_0$ qui
n'est pas de ramification,
alors il existe des coordonn\'ees locales $(x,y,z)$ au voisinage de
$\pi(R)$ telles que la surface
$S_0$ ait pour \'equation locale $xy=0$. L'\'eclatement $\pi$ \'etant
d\'efini par
$(x,w,z)\mapsto (x,y=xw,z)$, on v\'erifie que les
morphismes\hfill\break\vs{-4} $T_{\mu}$, $\gamma$
et
$\alpha$ ont respectivement pour matrices ${}^t[{1\atop 0}{0\atop 0}{0\atop
1}]$,
$[0,1,0]$ et $\biggl[{{\sc1\atop \sc0 }\atop\sc 0}{{\sc0\atop\sc x
}\atop\sc 0}{{\sc0\atop \sc0
}\atop
\sc1}\biggr]$. Pour la m\^eme raison que ci-dessus, on voit que le
morphisme $\zeta$ a pour matrice
$[x]$; ainsi la suite $(3.4.4))$ localis\'ee en $R$ est exacte, car
$\NN\hs{-1}\mid_{\lt_0}$ est
isomorphe
\`a $k[z]$, et $\lambda$ \'etant la r\'eduction modulo $x$.\cqfd\vs{2}
\sec{Corollaire}{On garde les notations de (3.4). Alors:\par
\ind 1) $\MM$ s'ins\`ere dans les suites exactes$$\leqalignno{
&0\to\MM\to\ost\to\cq\to0 &(3.5.1)\cr
&0\to\nstpt\to\MM\to\OO_{\lt_0}(2H-\lt_0)\to0& (3.5.2)\cr}$$
\ind 2)  $h^1({\MM})=0$, donc en particulier $h^0({\MM})=25$.}
\dem \par
\ind (1) Sachant que $\nstpt\simeq\OO_{{\st_0}}(4H-2\lt_0)$, alors la suite
$(3.5.1)$ est obtenue en
combinant les conclusions $(1)$ \`a $(3)$ de la proposition (3.4). La suite
$(3.5.2)$ est d\'eduite
directement du digramme (3.4.3).\par
\ind(2) La suite de cohomologie associ\'ee \`a $(3.5.2)$ donne en
particulier la suite exacte
$$0\to\cohr{0}{\nstpt}\to\cohr{0}{\MM}\to\cohr{0}{\OO_{\lt_0}(2H-\lt_0)}\to\
cohr{1}{\nstpt}
\to\cohr{1}{\MM}\to\cohr{1}{\OO_{\lt_0}(2H-\lt_0)}.$$  Or
$\cohr{1}{\OO_{\lt_0}(4H-\lt_0)}=0$ ( car
$4H-\lt_0$ est un diviseur de degr\'e $8$ sur une courbe elliptique), et
$h^1(\OO_{{\st_0}}(4H-2\lt_0))=h^1(\nstpt)=0$ ([W], 2.2), on en d\'eduit
donc
$h^1({\MM})=0$. Avec ceci, et en utilisant la suite de cohomologie
associ\'ee \`a $(3.5.1)$, on
obtient $h^0(\MM)\hs{-1}=h^0(\ost)-4$. Comme
$h^2(\ost)=h^0(\OO_{{\st_0}}(-4H))$ (par la dualit\'e de Serre), et
$h^0(\OO_{{\st_0}}(-4H))=0$ (car
$-4H.\Delta=-16<0$), alors, en utilisant le th\'eor\`eme de Riemann-Roch,
on a $h^0(\ost)={\dis
1\over\dis 2}(4H-\lt_0)(4H)+1 = 8 H^2 -2 \lt_0.H + 1=29$. On a donc
$h^0({\MM})=25$. \vs{3}
\sec{Proposition}{ Soient $S_0$, $L_0$, ${\st_0}$, $\lt_0$, $H$,
$R_1,\dots,R_4$ comme dans (3.2.3),
et $\MM$ le faisceau sur ${\st_0}$, d\'efini dans (3.4). Soit $C\subset
S_0$ une courbe lisse et
connexe, de degr\'e $d$ et genre
$g$, distincte de la droite $L_0$; on note $\ct$ sa transform\'ee
stricte.\par
\ind {\rm (i)} On a la suite exacte suivante:\par
 $\hs{40} 
0\to\cohr{0}{\MM(-\ct)}\to\cohr{0}{\MM}\to\cohr{0}{\nstct}\to\cohr{1}{\MM(-\
ct)}\to0$\par
\ind{\rm (ii)} Si $h^0(\OO_{{\st_0}}(4H-\ct-\lt_0))=0$, alors {\rm (a)}
$h^0(\MM(-\ct))=0$, et {\rm
(b)} il existe une injection $h^0(\cq)\hook\cohr{1}{\MM(-\ct)}.$\par
\ind {\rm (iii)} Si $h^1(\nctst)=0$, alors on a la suite exacte:\par 
$\hs{40}0\to\cohr{0}{\nctst}\to\cohr{0}{\ncp\hs{-1}}\to\cohr{0}{\nstct}\to
0.$}
\dem \par
\ind (i) La suite exacte
$0\to\OO_{{\st_0}}(-\ct)\to\OO_{{\st_0}}\to\OO_{\ct}\to 0$ tensoris\'ee par
le faisceau ( sans torsion) $\MM$ donne la suite exacte
$$0\to\MM(-\ct)\to\MM\to\nstct\to 0$$ d'o\`u l'on d\'eduit la suite
cherch\'ee dans (i) sachant que
$H^1(\MM)=0$ (cf. 3.5-(2)).\par
\ind (ii) On suppose $h^0(\OO_{{\st_0}}(4H-\ct-\lt_0))=0$. La suite (3.5.1)
tensoris\'ee par
$\OO_{{\st_0}}(-\ct)$ donne la suite exacte 
$$0\to H^0({\MM(-\ct)})\to H^0(\OO_{{\st_0}}(4H-\ct-\lt_0))\to H^0(\cq)\to
H^1({\MM(-\ct)})$$ d'o\`u
les assertions (a) et (b).\par
\ind (iii)  On a $C\simeq\ct$, car $C\not=L$, donc $\ncp\simeq\nctp\ $ o\`u
$\nctp=({\fitp\hskip -1mm\mid _{\ct}})/_{\dis\TT_{\ct}}$. Par ailleurs on a
$\nctst=(\TT_{{\st_0}}\mid_{\ct})/_{\dis \TT_{\ct}}$ et ${\MM\hskip
-1mm\mid_{\ct}}=(\fitp\mid_{\ct})/_{\dis (\TT_{{\st_0}}\mid_{\ct})}$. Ainsi
on a la suite exacte
$$0\to\nctst\to\nctp\to{\MM\hskip -1mm\mid _{\ct}}\fh0$$ L'hypoth\`ese $
h^1(\nctst)=0$ donne alors
la suite voulue.\cqfd\vfill\eject
{\bf Preuve du th\'eor\`eme 3.3 }\par
 Soit $C_0\subset S_0$ une courbe lisse et connexe v\'erifiant les
hypoth\`eses de (3.3). Il s'agit
de montrer l'in\'egalit\'e  $t_{(L_0,C_0,S_0)}\WW<t_{C_0}(H_{d,g})$. On
consid\`ere le morphisme
canonique
$f\colon \WW\to\QQ$ qui au triplet
$(L,C,S)$ associe la paire $(L,S)$. Soit $F$ la fibre de $f$ en
$(L_0,S_0)$, on a la suite exacte
des espaces tangents: $$0\to T_{(L_0,C_0,S_0)}F\to T_{(L_0,C_0,S_0)}\WW\to
T_{(L_0,S_0)}\QQ$$ La
fibre $F$ n'est autre que le sch\'ema de Hilbert $\Hdg(S_0)$ des courbes de
$\Hdg$ trac\'ees sur
$S_0$, donc $T_{(L_0,C_0,S_0)}F=\cohr{0}{\NN_{C_0/S_0}}$, et puisque 
$\dim_{C_0}
\Hdg(S_0)\leq \dim_{\ct_0}\Hdg({\st_0})$ (cf. [K1], lemme 22), on a
l'in\'egalit\'e
$$t_{(L_0,C_0,S_0)}\WW\leq
t_{(L_0,S_0)}{\QQ}+\dcoh{0}{\ct_0}{\NN_{\ct_0/\st_0}}.$$ Soit $\MM$ le
faisceau sur $\st_0$ introduit dans (3.4). D'apr\`es (3.5-(2)) on a
$h^0(\MM)=25$, et d'apr\`es le lemme (3.7) ci-dessous on a
$t_{(L_0,S_0)}{\QQ}=25$, donc
l'in\'egalit\'e pr\'ec\'edente devient
$$t_{(L_0,C_0,S_0)}\WW\leq h^0(\MM)+\dcohr{0}{\NN_{\ct_0/\st_0}}.$$ Or, on
a, d'une part
$\dcohr{0}{\NN_{\ct_0/\st_0}}=h^0(\NN_{C_0/\P{}{3}})-h^0(\MM\!\mid_{\ct_0})$
( cf. 3.6-(iii)), et
d'autre part $ h^0(\MM)\!=\!h^0(\MM\!\!\mid_{\ct_0})-h^1(\MM(-\ct_0))$ (cf.
3.6-(i),(ii)). Donc on
obtient, sachant que
$t_{C_0}(\Hdg)\!=\!h^0(\NN_{C_0/\P{}{3}})$, l'in\'egalit\'e 
$$t_{(L_0,C_0,S_0)}\WW\leq t_{C_0}(\Hdg)-h^1(\MM(-\ct_0))$$  Comme $
h^1(\MM(-\ct_0))\geq 4$ (cf.
3.6-(ii)), on a $t_{(L_0,C_0,S_0)}\WW< t_{C_0}(\Hdg)$, d'o\`u le
th\'eor\`eme.\cqfd\vs{3}
\sec{Lemme}{Avec les notations de (3.2.1) et (3.2.2), pour tout $(L,S)\in
\QQ$, on a $\dim
T_{(L,S)}\QQ=25$.} 
\dem Consid\'erons la suite exacte des espaces tangents
 $$0\to T_{S}(\QQ_L)\to T_{(L,S)}\QQ\to T_{L}(H_{1,0})\to0$$ induite par la
projection $\QQ\times
H_{1,0}\to H_{1,0}$. On a alors
$t_{(L,S)}\QQ=4+t_S(\QQ_L)$, car $H_{1,0}$ (la grassmannienne des droites
de $\P{}{3}$) est lisse de
dimension $4$. Comme, par d\'efinition,
$\QQ_L=\P{}{}\cohr{0}{\JJ_{L/\P{}{3}}^2(4)}$, il suffit donc de prouver que
$h^0(\JJ_{L/\P{}{3}}^2(4))=22$. Pour cela on calcule
$h^0(\JJ_{L^2/\P{}{3}}(4))$ o\`u $L^2$ est le
premier voisinage infinit\'esimal de $L$. Comme $L^2$ est li\'e \`a $L$ par
une intersection
compl\`ete de type $(1,2)$, il est arithm\'etiquement Cohen-Macaulay.
Consid\'erons alors la suite
exacte
$$0\fh\cohr{0}{\JJ_{L^2/\P{}{3}}(4)}\fh\cohr{0}{\JJ_{L/\P{}{3}}(4)}\fh\cohr{
0}
{\NN_{L/\P{}{3}}^{\vee}(4)}\to0$$  d'o\`u l'on tire
$h^0(\JJ_{L/\P{}{3}}^2(4))=h^0(\JJ_{L^2/\P{}{3}}(4))=h^0(\JJ_{L/\P{}{3}}(4))
-
h^0(\NN_{L/\P{}{3}}^{\vee}(4))$. On conclut sachant que
$h^0(\JJ_{L/\P{}{3}}(4))=30$ et que
$\dcohr{0}{\NN_{L/\P{}{3}}^{\vee}(4)}=\dcohr{0}{\OO_{L}(3)\oplus\OO_{L}(3)}=
8$.\cqfd\vs{3}
\sec{Remarques}{\rm Il y a une autre mani\`ere de prouver le th\'eor\`eme
(3.3) (cf. [Az], II-6): on
adapte
\`a notre situation les m\'ethodes de ([K2], I-1.3) utilisant la
cohomologie d'Andr\'e et la
th\'eorie d\'evelopp\'ee dans ([L]), pour calculer $T_{(L_0,C_0,S_0)}\WW$.
On montre
\smallskip
\hs{30}$T_{(L_0,C_0,S_0)}\WW=\fibp{T_{C_0}H_{d,g}}{\coh{0}{\ct_0}{\MM\mid_{\
ct_0}}}{\coh{0}{\st_0}{\MM}}$.
\smallskip
 Ensuite on consid\`ere la projection
$T_{(L_0,C_0,S_0)}\WW\to{T_{C_0}H_{d,g}}$, on montre que son
noyau (resp. conoyau ) est \'egal \`a $\cohr{0}{\MM(-\ct_0)}$ (resp. 
$\cohr{1}{\MM(-\ct_0)}$). Sachant que  
$\dcohr{0}{\MM(-\ct_0)}=0$ et $\dcohr{1}{\MM(-\ct_0)}\geq 1$, on
d\'eduit $t_{(L_0,C_0,S_0)}\WW<t_{C_0}(H_{d,g})$ comme voulu. }\vs{5}

\ind Dans ce qui suit, nous allons donner des conditions suffisantes afin
qu'une composante
irr\'eductible de $\Hdg$, contenant une courbe trac\'ee sur une quartique
\`a droite, soit
 domin\'ee par une composnte irr\'eductible de $\WW_{d,g}$. \par
\sec{Proposition}{Soit $S_0\subset\P{}{3}$ une surface quartique \`a droite
double $L_0$ comme dans
(3.2.3). Soit $C_0\subset S_0$ une courbe lisse et connexe, de degr\'e $d$
et genre $g$, qui coupe la
doite $L_0$ en $r$ points distincts des points-pinces et ce
transversalement. On note
$\ct_0$ la transform\'ee stricte de $C_0$ dans $\st_0$. Si\smallskip
 (3.9.1)\hs{30}$d>16$,\ \ $G(d,5)<g\leq (d-1)^2/8$,\ \ $r>36$,\smallskip 
ou bien si \smallskip
(3.9.2)\hs{15}$\cases{d>28,\ \  G(d,8)<g\leq (d-1)^2/8,\ \  r>72\cr {}\cr
\dcohr{0}{\OO_{\st_0}(-\w C_0)(7)}=\dcohr{1}{\OO_{\st_0}(-\w
C_0)(4)}=0.\cr}$
\smallskip
 Alors toute g\'en\'erisation de $C_0$ dans $\Hdg$ est trac\'ee sur une
surface quartique \`a droite
double.}
\dem Pour montrer la proposition, il revient au m\^eme de montrer que si
${\ee C}\subset\P{A}{3}$ est une d\'eformation plate de $C_0$, avec $A$ un
anneau de valuation
discr\`ete sur $k$, alors la  fibre g\'en\'erique de ${\ee C}$ au-dessus de
$A$ est une courbe
trac\'ee sur une quartique \`a droite double (i.e., le morphisme
classifiant $f\colon\spec A\to\Hdg$ 
se factorise via l'immersion $\VV\hook \Hdg$).\par 
\ind  Montrons d'abord que l'on peut se ramener au cas o\`u $A$ est
complet, i.e, $A$ est isomorphe
\`a un anneau de s\'eries formelles
$k[[t_1,\cdots,t_n]]$. Supposons qu'il existe un morphisme
$g\colon\spec
\hat A\to\VV$ , o\`u
$\hat A$ est le compl\'et\'e de $A$, faisant commuter le diagramme\vs{2}
$$\diagram{\spec \hat
A&\hfl{g}{}&\VV\cr
 \vfl{}{j}&&\vfl{}{}\cr
\spec A&\hfl{f}{}&\Hdg\cr}$$ Localement $\VV$ est le lieu des z\'eros d'un
id\'eal $(F_1,\dots
F_m)$, et la commutativit\'e du digramme ci-dessus
\'equivaut \`a dire que les $j^*f^* (F_n)\in \hat A$ sont nuls, et donc,
comme $j^*$ est une
injection, que les $f^* (F_n)$ sont nuls dans $A$; autrement dit le
morphisme $f$ se factorise via
$\VV$. Avec le m\^eme argument on peut montrer que l'on peut se ramener au
cas $A=k[[t]]$.\par
\ind Soit alors ${\ee C}\subset\P{A}{3}$ une d\'eformation plate de $C_0$,
avec $A=k[[t]]$. On
applique la proposition (2.9) (resp. (2.11)) si $C_0$ v\'erifie les
hypoth\`eses (3.9.1) (resp.
(3.9.2)), en prenant $s=4$ (resp. $s=4$ et $\sigma=7$). On en d\'eduit
l'existence d'une
d\'eformation plate
${\ee S}\subset\P{A}{3}$ de
$S_0$ contenant
$\ee C$ et dont la fibre g\'en\'erique est une surface singuli\`ere \`a
courbe double, qui est
g\'en\'erisation plate de
$L_0$, donc par platitude, c'est une droite double. \cqfd
\sec{Corollaire}{ Avec les notations de la proposition ci-dessus on a: 1)
$C_0$ est un point
int\'erieur de l'image de $\WW$ dans $\Hdg$. 2) Toute composante
irr\'eductible de $\Hdg$ contenant
$C_0$ est domin\'ee par une composante irr\'eductible de $\WW$.}
\dem Il suffit d'appliquer des r\'esultats plus g\'en\'eraux, par exemple
([DG], 1.10.2 et 1.10.3),
au morphisme
$\WW_{d,g}\to\Hdg$.\cqfd\vs{4}
\ind Comme con\'equence de tout ce qui pr\'ec\`ede on a le\par 
\sec{Th\'eor\`eme}{Soit $S_0$ une surface quartique \`a droite double $L_0$
comme dans (3.2.3).
Soient
$C_0$, $\ct_0$ comme dans la proposition (3.9). On suppose que  $C_0$
v\'erifie les hypoth\`eses
(3.9.1) ou (3.9.2). Si de plus
$\dcohr{1}{\NN_{\ct_0/\st_0}}=\dcohr{0}{\OO_{\st_0}(4H-\lt_0-\ct_0)}=0$, 
alors toute composante
irr\'eductible de $\Hdg$ contenant $C_0$ est non r\'eduite.}
\dem Si $C_0$ v\'erifie (3.9.1) ou (3.9.2), alors d'apr\`es (3.10), toute
composante irr\'eductible
de
$\Hdg$ contenant $C_0$ est domin\'ee par une composante irr\'eductible de
$\WW_{d,g}$, c'est donc
une composante non r\'eduite d'apr\`es (3.3).\cqfd\vs{3} 
\centerline{\bf Preuve du th\'eor\`eme (3.1) }
\ind Nous arrivons maintenant \`a la preuve du th\'eor\`eme (3.1). Elle
consiste \`a montrer que pour
tout couple d'entiers $(d,g)$ satisfaisant aux hypoth\`eses dudit th\'eor\`eme, on peut construire
une courbe lisse et connexe de degr\'e  $d$ et genre $g$ trac\'ee sur la
surface $S_0$ satisfaisant
aux conditions (3.9.2) de la propopsition (3.9). Ceci sera facilit\'e par
le fait qu'il y a des
moyens num\'eriques pour classifier les courbes de $\st_0$, et donc de
transformer les conditions
(3.9.2) en termes purement arithm\'etiques. L'essentiel de la preuve sera
de rendre effectives ces
conditions num\'eriques. Pour ce faire, on commencera d'abord par prouver
la
\par
\sec{Proposition}{Soient $S_0$, $\st_0$, $\Delta$, $E_1,\cdots,E_9$, $L_0$,
$\lt_0$, $H$ comme dans
(3.2.3). Soit
$D=\delta\Delta-\sum_{i=1}^9 m_iE_i$ un diviseur de $\st_0$ tel que \par
(1) $\delta\geq
m_1+m_2+m_3$\par  (2) $ m_1>m_2\geq m_3\geq \dots  \geq m_9\geq1 $\par (3)
$r:=3\delta -
\sum_{i=1}^9 m_i\geq 3$\par Alors $D$ est effectif et le syst\`eme
lin\'eaire $\mid\!\!D\!\!\mid$
contient une courbe lisse et connexe $\ct$ telle que \par  (a)
$H^1(\NN_{\ct/\st_0})=0$ et
$\lt_0\!\cap\!\ct$ ne contient pas les points de ramification ni aucune
paire de points en
involution; le cardinal de $\ct\cap\lt_0$ est \'egal \`a $r$. \par (b) Si
$\delta>13$ alors
$\cohr{0}{\OO_{\st_0}(4H-\lt_0-\ct)}=0$.\par (c) Si $\delta>28$ alors
$\cohr{0}{\OO_{\st_0}(7H-\ct)}=0$.\par (d) Si $m_9\geq 4$, $m_1\geq m_2+4$,
$r\geq 9$ et
$\delta>m_1+10$, alors
$\cohr{1}{\OO_{\st_0}(4H-\ct)}=0$.\par
\ind Soit $C$ l'image par $\varphi$ de $\ct$, alors est une courbe lisse et
connexe, isomorphe \`a
$\ct$, coupant la droite $\L_0$, transversalement, en $r$ points distincts
des points de
ramification. Le degr\'e
$d$ de $C$ est donn\'e par $d=4\delta-2m_1-\sum_{i=2}^9 m_i$ et son genre
$g$ par
$g=[(\delta-1)(\delta-2)-\sum_{i=1}^9 m_i(m_i-1)]/2$  }
\dem En appliquant ([bH], 3.1 et 3.4), les conditions conditions (1) \`a
(3) impliquent que
$\mid\!\!D\!\!\mid$ est sans points bases, donc (compte tenu du
th\'eor\`eme de Bertini) un
\'element g\'en\'eral
$\ct$ de $\mid\!\!D\!\!\mid$ est une courbe lisse qui \'evite les points de
ramification et, en
raison de la condition (3) (qui \'equivaut \`a $\ct.\lt_0\geq 3$), on peut
voir qu'elle \'evite
aussi les points en involution. En vertu de ([bH], 1.1),
$\cohr{1}{\OO_{\st_0}(-\ct)}=0$,
donc, en fait,
$\ct$ est connexe. Le nombre de points de $\ct\cap\lt_0$ est \'egal \`a
$\ct.\lt_0=3\delta -
\sum_{i=1}^9 m_i=r$. \par
\ind Soit $g$ le genre de $\ct$; la formule d'adjonction appliqu\'ee \`a
$\ct$ et la condition (3)
permettent de voir que  deg $\NN_{\ct/\st_0}$$=\ct^2>2g-2$ et donc
$\cohr{1}{\NN_{\ct/\st_0}}=0$.\par
\ind Pour voir (a), on observe que $\delta>13$ implique $\Delta.\ct<0$,
donc, d'apr\`es ([D], p.24),
$\cohr{0}{\OO_{\st_0}(4H-\lt_0-\ct)}=0$. Le m\^eme argument permet de voir
que
$\cohr{0}{\OO_{\st}(4H-\lt-\ct)}=0$ si $\delta>28$, d'o\`u l'assertion
(c).\par
\ind Reste \`a montrer (d). Pour cela on rappelle le r\'esultat suivant
d\^u \`a M. Skiti, (cf.
[S]):\par
\ind {\bf Lemme.}{\it  Soit $\LL$ un fibr\'e en droites sur $\st_0$
correspondant, dans la base
$(\Delta, -E_1,\dots,-E_9)$, au multi-entier $(a, b_1,\dots ,b_9)$ avec
$a>b_1\geq b_2\geq\cdots\geq
b_9$, $a\geq b_1+ b_2+b_3$ et $\LL.\lt_0>0$. Si $b_9\geq -1$ alors
$\cohr{1}{\LL}=0$.}\vs{2}Posons
$\LL= \OO_{\st_0}(-\lt_0-4H+\ct)$. Le multi-entier entier associ\'e \`a
$\LL$ est \'egal
\`a $ (\delta- 19,m_1-9,m_2-5,\dots  ,m_9-5)$. Ainsi sous les hypoth\`eses
de (d), $ (\delta-
19,m_1-9,m_2-5,\dots  ,m_9-5)$  v\'erifie celles du lemme ci-dessus, donc 
$\cohr{1}{\OO_{\st_0}(\LL)}=0$, et par la dualit\'e de Serre
$\dcohr{1}{\OO_{\st_0}(\LL)}=\dcohr{1}{\OO_{\st_0}(4H-\ct)}$, d'o\`u
l'assertion (d).\par
\ind La derni\`ere partie est assez claire: le calcul de $d$ et $g$ est
classique, et le cardinal de
$C\cap L_0$ \'etant \'egal \`a $\ct.\lt_0$ vaut $r$; les autres
affirmations d\'ecoulent de
([H3], 3.1).\cqfd\vs{5}
\sec{Preuve du th\'eor\`eme (3.1) }{}
\ind Soient $d$ et $g$ deux entiers satisfaisant aux hypoth\`eses $(A)$ et
$(B)$, alors d'apr\`es
le corollaire (3.16) ci-dessous, il existe un existe un d\'ecuple d'entiers
$(\delta,m_1,\cdots,m_9)$  v\'erifiant les  in\'egalit\'es:\par 
\hs{10}{\rm (i)}\  \  $m_1\geq m_2\geq\cdots\geq m_9\geq 4$\par
\hs{10}{\rm (ii)}\  \  $\delta\geq m_1+m_2+m_3 $\par
\hs{10}{\rm (iii)}\  \  $m_1\geq m_2+4$\smallskip
\hs{10}{\rm (iv)}\  \  $3\delta-\sum_{i=1}^{9}m_i\geq 73$\par tel que $d$
et $g$ soient donn\'es par
les formules\par
  $d=4\delta-2m_1-\sum_{i=2}^{9}m_i$, 
 $g=[(\delta-1)(\delta-2)-\sum_{i=1}^{9}m_i(m_i-1)]/2$.\par
On consid\`ere le diviseur $D=\delta\Delta-\sum_{i=1}^9 m_iE_i$, alors $D$
v\'erifie, en
particulier, les hypoth\`eses (1), (2), (3), (b), (c) et (d) de la
proposition (3.12). Ainsi il
existe une courbe lisse et connexe $\ct_0$ dans le syst\`eme lin\'eaire
$\mid\! D\!\mid$ de
genre $g$ satisfaisant aux hypoth\`eses du th\'eor\`eme (3.11). Et ceci
termine la d\'emonstration
du  th\'eor\`eme.\cqfd\vs{3}
\centerline{\bf Appendice arithm\'etique}\vs{2}
Nous \'etablissons dans ce qui suit les arguments arithm\'etiques
utilis\'es dans la preuve du
th\'eor\`eme (3.1). Les m\'ethodes de d\'emonstration sont largement
inspir\'ees de celles
utilis\'ees par Gruson et Peskine dans ([GP2]) et reprises par Hartshorne
dans ([H3]).\vs{2}

\sec{Lemme}{ Soit $d$ un entier $\geq 95$. Pour tout entier $g\in\rbrack\
G(d,8)\ ,{(d-1)^2/
8}\rbrack$ il existe un entier $v \in\lbrack d/8+9,  (d+1)/2\rbrack$, et
des demi-entiers, (i.e. des
\'el\'ements de 
${\dis{1\over2}{\bb Z}}$), $\alpha_2,\dots,\alpha_9$ v\'erifiant:\smallskip

\hs{10}{\rm (1)}\  \    $F_d(v-1)<g\leq F_d(v)$\ , \ avec
$F_d: x\mapsto{\dis {1\over 2}}(x-1)(d-x)$\ \par
\hs{10}{\rm (2)}\  \   
$\mid\alpha_2\mid\leq\alpha_3\leq\cdots\leq\alpha_9\leq {\dis {v\over
2}}-4$\par
\hs{10}{\rm (3)}\  \   $-\alpha_2+\alpha_3+\cdots+\alpha_9\leq d-v-8$\par
\hs{10}{\rm (4)}\  \    $2\alpha_i\equiv v\mod{2}$\par
\hs{10}{\rm (5)}\  \    $d-\sum_{i=2}^{9}\alpha_i \equiv 0\mod{2}$\par tels
que\par \hs{10} {\rm
(6)}\ \ $g=F_d(v)+1-{\dis {1\over 2}}\sum_{i=2}^{9}\alpha_i^2$.}
\sec{Remarque}{\rm  Il est facile de voir que les congruences (4) et (5)
impliquent que le nombre $
F_d(v)-{\dis {1\over 2}}\sum_{i=2}^{9}\alpha_i^2$ est bien un entier. Par
ailleurs
 la fonction $F_d$ atteint son maximum pour $x=(d+1)/2$, et vaut dans ce
cas $(d-1)^2/8$, donc la
formule (5) n'a de sens que si l'on impose $g\leq (d-1)^2/8$, ce qui est le
cas ici.} 
 {\it Preuve du lemme. } Soit $d$ un entier $\geq 95$. Soit $g$ un entier
dans $\rbrack\
G(d,8),(d-1)^2/8\rbrack$. Il est facile de voir que $F_d(d/8+9)\leq
G(d,8)$, et puisque la fonction
$F_d$ est bijective sur l'intervalle $[d(d/8+9),(d+1)/2]$, alors ce dernier
contient un r\'eel
$\beta$ tel que $F_d(\beta)=g$. Si $\beta\in \bb N$, on pose $v=\beta$. Si
$\beta\not\in \bb N$, on a
$[\beta]+1\leq (d+1)/2$, ($[x]$ symbolise la  partie enti\`ere de $x$); on
pose  $v=[\beta]+1$. Il
est clair que $v$ v\'erifie la double in\'egalit\'e (1), compte tenu de la
stricte croissance de
$F_d$ dans $[d(d/8+9),(d+1)/2]$.\par Posons $n=2(F_d(v)+1-g)$. Alors $n$
est un entier, et, en
notant que
$2(F_d(v)-F_d(v-1))=d-2v+2$, on a $2\leq n<d-2v+4$ et donc  
 $2\leq n<6v-68$. On va montrer qu'un tel $n$ s'\'ecrit
$n=\sum_{i=2}^{9}\alpha_i^2$, o\`u les
$\alpha_i$ sont dans
$(1/2)\bb Z$ v\'erifiant les conditions (2) \`a (5) du lemme. Pour ce
faire, on va distinguer deux
cas, selon la parit\'e de
$v$.\par
\ind Supposons $v$ pair. Il s'agit de trouver 8 entiers $\alpha_i$ avec
$\mid\alpha_i\mid\leq v/2-4$
tels que $n=\sum_{i=2}^{9}\alpha_i^2$. Posons $q=v/2-4$, alors
$n<3q^2-2q+3$. Ceci implique,
d'apr\`es ([H3], 2.1), que
$n$ est somme de 5 carr\'es d'entiers, $\alpha_5,\cdots,\alpha_9$ avec
$\mid\alpha_i\mid\leq N$, et
donc, a fortiori, somme de 8 carr\'es d'entiers en posant
$\alpha_2=\alpha_3=\alpha_4=0$.\par
\ind  Il reste \`a verifier la congruence (5) et l'in\'egalit\'e (3). (
L'in\'egalit\'e (2) \'etant
\'evidemment satisfaite, quitte \`a r\'earranger les $\alpha_i$.). Puisque
$g$ est un entier, on a
$(d-v)(v-1)\equiv \sum_{i=2}^{9}\alpha_i^2\mod{2}$, ou encore
$d\equiv\sum_{i=2}^{9}\alpha_i^2\mod{2}$, d'o\`u la
congruence (5). L'in\'egalit\'e (3) est \'evidente car $v\geq12$ et $\sum_{i=2}^{9}\alpha_i^2=n\leq
d-2v+4$.\smallskip
\ind Si maintenant $v$ est impair, on a $F_d(v)\in\bb Z$. Alors $m=4n-3$
satisfait les hypoth\`eses
de ([H3], 2.2),i.e., $m\equiv 5\mod{8}$ et $m<3p^2+2p+1$ avec $p=v-8$.
D'o\`u l'existence de $5$
entiers impairs $\beta_i$ avec $\mid\beta_i\mid\leq v-8$ tels que
$m=\sum_{i=5}^{9}\beta_i^2$. On pose
$\mid\alpha_2\mid=\alpha_3=\alpha_4=1/2$ et
$\alpha_i=\mid\beta_i/2\mid$ pour $i\geq 5$. Le choix du signe de
$\alpha_2$ se fait ainsi: on
\'ecrit $\sum_{i=5}^{9}\beta_i=2M+1$, ensuite on prend $\alpha_2=1/2$
(resp. -1/2) si $M\equiv d\mod{2}$ (resp.
$M\equiv d-1\mod{2}$), et ceci permet d'obtenir la congruence (5). Il reste
\`a v\'erifier
l'in\'egalit\'e  (3); pour cela il suffit de v\'erifier
$\sum_{i=5}^{9}\beta_i\leq 2d-2v-19$. Puisque 
$m=\sum_{i=5}^{9}\beta_i^2$ et $m+5<4d-8v+18<2(2d-2v-19)$ (car $v\geq 14$),
il suffit de v\'erifier
$\sum_{i=5}^{9}\beta_i\leq (\sum_{i=5}^{9}\beta_i^2)+5/2$, mais ceci n'est
autre que la lapalissade
$\sum_{i=5}^{9}(\beta_i-1)^2\geq 0$. Ceci d\'emontre le lemme.\cqfd
\ind Si dans l'\'enonc\'e du lemme ci-dessus on remplace $v\leq (d+1)/2$
par $v\leq d-73$, (comme
nous serons amen\'e \`a le faire en vue de montrer (3.16)), alors celui-ci
reste valable pour $d\geq
147$ puisque, dans ce cas, $(d+1)/2\leq d-73$, et pour qu'il le reste
lorsque $95\leq d\leq 146$, (
auquel cas $d-73<(d+1)/2$), il faut \'evidemment remplacer dans la preuve
pr\'ec\'edente
l'hypoth\`ese
$g\leq(d-1)^2/8$ par $g\leq F_d(d-3)=73(d-74)/2$.  Ceci explique la
pr\'esence des r\'egions $(A)$
et $(B)$ dans le prochain corollaire. L'int\'eret de la condition $d\geq95$
s'explique par le fait
que si l'on cherche un entier naturel
$v\leq d-73$, il faut a priori avoir $d\geq 73$, (donc on conna\^it
$G(d,8)$), et on montre moyennant
un calcul \'el\'ementaire que $F_d(d-73)<G(d,8)$ pour $d\leq 94$. Ainsi,
puisque dans ce cas la
fonction
$F_d(x)$ est strictement croissante pour $x\leq d-73$ , l'\'egalit\'e (6)
n'est jamais r\'ealisable  
avec $d\leq 94$ et $g>G(d,8)$.\vs{4}   
\sec{Corollaire}{Soient $d$ et  $g$ deux entiers satisfaisant:\smallskip
$(A)$\hs{30}\ \ $G(d,8)<g\leq {\dis {73\over 2}}(d-74) $\ \  si\ \ $146\geq
d\geq 95$\smallskip
$(B)$\hs{30} \ \  $G(d,8)<g\leq {\dis{(d-1)^2\over 8}}$\ \  \ si\ \
$d\geq147$\smallskip Alors il
existe un d\'ecuple d'entiers $(\delta,m_1,\cdots,m_9)$  v\'erifiant les 
in\'egalit\'es:\par 
\hs{10}{\rm (i)}\  \  $m_1\geq m_2\geq\cdots\geq m_9\geq 4$\par
\hs{10}{\rm (ii)}\  \  $\delta\geq m_1+m_2+m_3 $\par
\hs{10}{\rm (iii)}\  \  $m_1\geq m_2+4$\smallskip
\hs{10}{\rm (iv)}\  \  $3\delta-\sum_{i=1}^{9}m_i\geq 73$\par tel que $d$
et $g$ soient donn\'es par
les formules\par
\hs{10}{\rm (v)}\  \  $d=4\delta-2m_1-\sum_{i=2}^{9}m_i$\par
\hs{10}{\rm (vi)}\  \  $2g=(\delta-1)(\delta-2)-\sum_{i=1}^{9}m_i(m_i-1)$.}

\dem Soient $d$  et $g $ deux entiers satisfaisant aux hypoth\`eses $(A)$
ou $(B)$ du corollaire. On
sait, d'apr\`es le lemme (3.15) et le commentaire qui le suit, qu'il existe
un entier $v \in\lbrack
d/8+9,\ d-73\rbrack\ $, et des demi-entiers $\alpha_2,\cdots,\alpha_9$
v\'erifiant les  conditions
(i) \`a (vi) du lemme. On pose  \smallskip
\hs{15}\ \ $\delta={\dis {1\over 2}} (d+2v-\sum_{i=2}^{9}\alpha_i)$\ ,\ \ 
$m_1=\delta-v$
\ ,\  \  $m_i={\dis {v\over 2}}-\alpha_i$\ \ pour $i\geq 2\ .$\smallskip
Alors les congruences (4) et
(5) du lemme impliquent respectivement $\delta\in\bb Z$ et
$m_i\in\bb Z$ pour $i\geq 2$, et une v\'erification simple permet de voir
que le multi-entier
$(\delta,m_1,\cdots,m_9)$ ainsi d\'efini satisfait aux conditions (1) \`a
(6), ce qui d\'emontre le
corollaire.\cqfd\vs{5}{\bf Bibliographie}\par

\ref{Ar}{M. Artin}{Deformations of singularities}{Tata Institute Lectures,
54, 1976}
\ref{Az}{S. Azziz}{Exemples de composantes irr\'eductibles non r\'eduites
du sch\'ema de Hilbert des
courbes lisses et connexes de $\P{}{3}$}{Th\`ese, universit\'e de
Toulouse-III, 1996}
\ref{D}{M. Demazure}{Surfaces de Del Pezzo}{Lect. Notes in Math. 777
(1980), 23-69}
\ref{DG}{J. Dieudonn\'e, A. Grothendieck}{Etude locale des sch\'emas et des morphismes de
sch\'emas}{EGA IV, Publ. I.H.E.S, 20 (1964) }
\ref{DP}{A. Dolcetti, G. Pareschi}{On linearly normal space
curves}{Math.Zeit. 198 (1988), 73-82}
\ref{E}{Ph. Ellia}{D'autres composantes non r\'eduites de Hilb
$\P{}{3}$}{Math. Ann. 277 (1987),
433-446} 
\ref{EF}{Ph. Ellia, M. Fiorentini}{D\'efaut de postulation et
singularit\'es du sch\'ema de
Hilbert}{Ann. Uni. Ferrara Nuova, Ser.Sez.VII, 30 (1984), 185-198}
\ref{GP1}{L. Gruson, Ch. Peskine}{Genre des courbes alg\'ebriques dans
l'espace projectif}{Lect.
Notes in Math. 687 (1978), 31-59}
\ref{GP2}{L. Gruson, Ch. Peskine}{Genre des courbes alg\'ebriques dans
l'espace projectif (II)}{Ann.
Scient. \'Ec. Norm. Sup. 15 (1982), 410-418}  
\ref {H1}{R. Hartshorne}{Une courbe irr\'eductible non lissifiable dans
$\P{}{3}$}{C.R. Acad. Sci. 
Paris, t. 299, s\'erie I, 5 (1984), 133-136}
\ref{H2}{R. Hartshorne}{Algebraic Geometry}{Gradu. Texts in Math.,
Springer-Verlag, 52 (1977)}
\ref{H3}{R. Hartshorne}{Genre des courbes alg\'ebriques dans l'espace
projectif}{S\'em. Bourbaki,
592 (1981/82), 1-13}
\ref {H4}{R. Hartshorne}{Stable Reflexive Sheaves}{Math. Ann. 254 (1980),
121-176}
\ref {bH}{B. Harbourne}{Complete Linear Systems On Rational
Surfaces}{Trans. A.M.S, 289, 1 (1985),
213-226}
\ref {K1}{J. Kleppe}{Non reduced components of the Hilbert scheme of smooth
space curves}{Lect.
Notes in Math.1266 (1987) 181-207}
\ref {K2}{J. Kleppe}{The Hilbert-Flag scheme... }{Th\`ese, universit\'e
d'Oslo (1982)}
\ref{Ka1}{N.M. Katz}{Pinceaux de Lefschetz: th\'eor\`eme d'existence}{Lect.
Notes in Math. 340,
(1973), expos\'e 17, 212-253}
\ref{Ka2}{N.M. Katz}{Etude cohomologique des pinceaux de Lefschetz}{Lect.
Notes in Math. 340,
(1973), expos\'e 18, 254-327}
\ref {L}{ O.A. Laudal}{Formal moduli of algebraic structures}{Lect. Notes
in
Math. 754 (1979)} 
\ref{M1}{D. Mumford}{Further pathologies in algebraic geometry}{Amer. J.
Math. 84 (1962), 642-648}
\ref{M2}{D. Mumford}{Red book of varieties and schemes}{Lect. Notes in
Math. 1358 (1988)}
\ref{R}{D.S. Rim}{Formal deformation theory}{Lect. Notes in
Math. 288 (1972), expos\'e 6, 32-132}
\ref{Ro}{J. Roberts}{Singularity subschemes and generic projections.
}{Trans. A.M.S, 212
(1975) 229-268}
\ref{S}{M. Skiti}{Courbes trac\'ees sur les surfaces quartiques \`a droite
double. }{Th\`ese,
universit\'e  Lille I, 1986}
\ref{T1}{B. Teissier}{Vari\'et\'es polaires I}{Invent. Math. 40 (1970),
267-292}
\ref{T2}{B. Teissier}{R\'esolution simultan\'ee: I- familles de
courbes}{S\'em. sur les
singularit\'es des surfaces, Ecole Polyth., Paris (D\'ecembre 1976)}
\ref{W}{Ch. Walter}{Curves on surfaces with multiple line}{J. reine angew.
Math. 412 (1990), 48-62}
\vs{15}\hs{100} Toulouse, Janvier 1997.
\vs{20}
Sa\"id AZZIZ\par
Laboratoire \'Emile Picard\par
Universit\'e Paul Sabatier,\par
118, Route de Narbonne. 31062 Toulouse. France.\par
e-mail: azziz@picard.ups-tlse.fr
\end